\documentclass[11pt]{article}
\usepackage{sectsty,latexsym,amssymb,dsfont,textcomp}
\usepackage[centertags]{amsmath}
\usepackage{fullpage}
\usepackage{color}
\usepackage{epsfig,psfrag}

\usepackage[centertags]{amsmath}
\sectionfont{\normalsize}
\subsectionfont{\normalsize}
\begin{document}
\date{}
\numberwithin{equation}{section}
\title{Singular limiting induced from continuum solutions \\
and the problem of dynamic cavitation
}

\author{
Jan Giesselmann\footnote{Weierstrass Institute, Berlin, Germany,
     jan.giesselmann@wias-berlin.de} 
%
\and 
Athanasios E. Tzavaras\footnote{Department of Applied Mathematics, University of 
Crete, Heraklion, Greece and Institute for Applied and Computational Mathematics, FORTH\@, Heraklion, Greece,
tzavaras@tem.uoc.gr}
}

\maketitle

\newtheorem{lemma}{Lemma}[section]
\newtheorem{theorem}[lemma]{Theorem}
\newtheorem{proposition}[lemma]{Proposition}
\newtheorem{maintheorem}[lemma]{Main Theorem}
\newtheorem{corollary}[lemma]{Corollary}
\newtheorem{definition}[lemma]{Definition}
\newtheorem{remark}[lemma]{Remark}
\newtheorem{remarks}[lemma]{Remarks}
\newtheorem{Notation}[lemma]{Notation}
\newcommand{\proof}{\noindent {\it Proof}\;\;\;}
\newcommand{\qed}{\protect~\protect\hfill $\Box$}

\newcommand{\wkarr}{\; \rightharpoonup \;}
\def\Weak{\,\,\relbar\joinrel\rightharpoonup\,\,}

\font\msym=msbm10
\def\Real{{\mathop{\hbox{\msym \char '122}}}}
\def\R{\Real}
\def\charf {\mbox{{\text 1}\kern-.24em {\text l}}} 

\def\A{\mathbb A}
\def\Z{\mathbb Z}
\def\K{\mathbb K}
\def\J{\mathbb J}
\def\L{\mathbb L}
\def\D{\mathbb D}
\def\cD{\mathcal D}
\def\cO{\mathcal O}
\def\cQ{\mathcal Q}
\def\Mink{{\mathop{\hbox{\msym \char '115}}}}
\def\Integers{{\mathop{\hbox{\msym \char '132}}}}
\def\Complex{{\mathop{\hbox{\msym\char'103}}}}
\def\C{\Complex}
\font\smallmsym=msbm7

\newcommand{\del}{\partial}
\newcommand{\eps}{\varepsilon}
\def\div{\hbox{div}\,}
\def\supp{\hbox{supp}\,}
\def\dist{\hbox{dist}\,}
\newcommand{\tcb}{\textcolor{blue}}
\newcommand{\tcg}{\textcolor{green}}
\newcommand{\tcr}{\textcolor{red}}

%
%

%

\begin{abstract}
\noindent
\noindent In the works of
K.A. Pericak-Spector and S. Spector \cite{ps88, ps98} a class of self-similar
solutions are constructed for the equations of radial isotropic elastodynamics
that describe cavitating solutions. Cavitating solutions decrease the total
mechanical energy and provide a striking example of non-uniqueness of entropy weak solutions
(for polyconvex energies) due to point-singularities at the cavity. To resolve 
this paradox, we introduce the concept of singular limiting induced from continuum solution (or slic-solution),
according to which a discontinuous motion is a slic-solution  if its averages 
form a family of smooth approximate solutions to the problem.  It turns out that there is an energetic cost for
creating the cavity, which is captured by the notion of slic-solution but neglected by the
usual entropic weak solutions.  Once this cost is accounted for, the total mechanical energy of the
cavitating solution is in fact larger than that of the homogeneously deformed state.
We also apply the notion of slic-solutions to a one-dimensional example describing the onset of fracture,
and to gas dynamics in Langrangean coordinates with Riemann data inducing vacuum in the wave fan.

\end{abstract}
%
%
%
\section{Introduction}
\label{intro}
\setcounter{equation}{0}

The goal of this study is to re-assess an example of non-uniqueness of entropy weak solutions for multi-dimensional
systems of hyperbolic conservation laws constructed by K.A. Pericak-Spector and S. Spector \cite{ps88, ps98}.
The example is associated with the the onset of cavitation -- from a homogeneously deformed state --
for the equations of nonlinear elasticity in homogeneous and isotropic elastic media.
Any attempt to study solutions that lie at the limits of continuum modeling (like cavities or shear bands)
needs to reckon with the problem of giving a proper definition for such solutions. Proposing such a definition 
is the main premise of the present work.

We consider the initial-boundary value problem 
for the equations of  elasticity 
\begin{align}
 \label{eq:introelas3d}
{\bf y}_{tt} - \operatorname{div} {\bf S} (\nabla {\bf y}) =0 , \;  & \\
\label{eq:intro-ic3d}
{\bf y}({\bf x},0 ) = \lambda {\bf x} &  \\
\label{eq:intro-bc3d}
{\bf y}({\bf x},t) = \lambda {\bf x} & \text{ for } |{\bf x}| > {\bar r} t  \, , 
\end{align}
for a given stretching  $\lambda>0$ and some ${\bar r}  >0.$
Here, ${\bf y}:\R^d \times \R_+ \rightarrow \R$ stands for the motion, ${\bf F} = \nabla {\bf y}$ is the
deformation gradient, and we 
employ the constitutive theory of hyperelasticity, 
that the Piola--Kirchhoff stress ${\bf S}$ is given as the gradient of a stored energy density
\begin{equation}
 \label{eq:energy3d1}
{\bf S} ({\bf F})= \frac{\partial W}{\partial {\bf F}} ({\bf F}) \, , 
\quad W: \R^{d \times d}_+ := \{{\bf F} \in \R^{d \times d} : \det({\bf F})>0\}  \longrightarrow \R \, .
\end{equation}
The homogeneous deformation $\bar{\bf y}({\bf x},t) = \lambda {\bf x}$ solves
\eqref{eq:introelas3d}-\eqref{eq:intro-bc3d}.

The stored energy function  determines the constitutive properties of the elastic material and - due to
frame indifference - has to be invariant under rotations. For homogeneous and {\it isotropic} elastic
materials $W$ takes the simplified form $W( {\bf F})=\Phi(\lambda_1, ... , \lambda_d)$, where $\Phi$ is a symmetric
function of the eigenvalues  $\lambda_1, ... , \lambda_d$  of $\sqrt{{\bf F \bf F}^ T}$; see \cite{TN65}. 
In that case \eqref{eq:introelas3d} admits as solutions  radially symmetric motions,
\begin{equation}
 \label{eq:introradsym}
{\bf y}({\bf x},t) = w(|{\bf x}|,t) \frac{\bf x}{|{\bf x}|}  \, , 
\end{equation}
generated by solving for the amplitude $w: \R_+ \times \R \longrightarrow \R_+$ the
initial-boundary value problem for the equations of isotropic radial elastodynamics,
\begin{align}
 \label{eq:weq}
w_{tt} &= \frac{1}{R^{d-1}} \partial_R  \left( R^{d-1} \frac{\partial \Phi }{\partial \lambda_1} (w_R,\frac{w}{R},\dots,\frac{w}{R})\right)
- \frac{1}{R} (d-1) \frac{\partial \Phi} {\partial \lambda_2} (w_R,\frac{w}{R},\dots,\frac{w}{R}).
\\
\label{eq:introwcond}
&\qquad \left \{
\begin{aligned}
w(R,0) &= \lambda R \, ,   \\
w(R,t) &= \lambda R \quad \text{ for } R > {\bar r} t \, .
\end{aligned}
\right .
\end{align}
This problem admits the special solution $\bar w(R,t) = \lambda R$ corresponding to a state of homogeneous deformation. 
The question arises if additional solutions of \eqref{eq:introelas3d}-\eqref{eq:intro-bc3d} may be constructed
by solving the problem \eqref{eq:weq}-\eqref{eq:introwcond}.

This idea was pursued by Ball \cite{ball82} using methods from the calculus of variations  to construct cavitating solutions
 in the context of static elasticity.
The study \cite{ball82} was groundbreaking in part because continuum modeling is used to address a problem that at
least {\it pro-forma} concerns situations beyond its range of applicability. 
There exists a critical stretching $\lambda_{cr}$ so that for $\lambda < \lambda_{cr}$ the only
equilibrium solution is the homogeneously deformed state; by contrast, for $\lambda > \lambda_{cr}$ there exist
non-trivial equilibria corresponding to a cavity in the material and with energy 
{\it less} than the energy of the homogeneous deformation  \cite{ball82}. 
 The reader is referred to 
 \cite{SS03,L09,NS11} (and references therein) for an account of developments concerning 
 cavitating equilibrium solutions in nonlinear elasticity.

In another remarkable development, K.A. Pericak-Spector and S. Spector \cite{ps88, ps98} use the ansatz
\begin{equation}
 \label{eq:intrordef}
w(R,t)= t \, r \Big ( \frac{R}{t} \Big ) \, ,
\end{equation}
to construct a self-similar weak solution for the dynamic problem \eqref{eq:introradsym}-\eqref{eq:introwcond}
that corresponds to a spherical cavity emerging at time $t=0$ from a homogeneously deformed state.
The solution in \cite{ps88} is constructed in dimension $d \ge 3$ for polyconvex stored 
energies of the special form
\begin{equation}
 \label{eq:energy3d2}
W({\bf F}) = \frac{1}{2} \sum_{i=1}^d \lambda_i^2 + h\left(\prod_{i=1}^d \lambda_i\right)
\tag {H$_1$}
\end{equation}
where $h: \R_+ \longrightarrow \R_+$ satisfies the hypotheses
\begin{equation}
 \label{eq:defh}
 h''>0 \, ,  \quad  h'''<0 \, \quad  \lim_{v \rightarrow 0} h(v) = \lim_{v \rightarrow \infty} h(v) = \infty  \, .
 \tag {H$_2$}
\end{equation}
The hypothesis $h'' >0$ refers to polyconvexity, while $h''' < 0$ indicates elasticity of softening type;
more general stored energies were treated in \cite{ps98}.
It is proved in \cite{ps88,ps98} that the self-similar solution has smaller mechanical energy than 
the associated homogeneously deformed state from where it emerges, and thus provides an example of
nonuniqueness of entropy weak solutions (at least for polyconvex energies). As already noted in \cite{ps88},
the paradox arises that by opening a cavity the energy of the material decreases, 
what induces an autocatalytic mechanism for failure.

There is a class of problems in material science, such as fracture, cavitation or shear bands,
where discontinuous motions emerge from smooth motions
via a mechanism of material instabilities.
 Of course after the material breaks or a shear band forms the motion can no longer be described at the level
of continuum modeling and microscopic modeling or higher-order regularizing mechanisms have to be taken into account.
Still, as such structures develop there is expected an intermediate time scale where both types of modeling apply.
The premise of this work is to introduce a meaning for solutions at this intermediate scale and to explore its ramifications.
The idea is roughly the following, presented here at the level of \eqref{eq:introelas3d}.
Given a possibly discontinuous motion ${\bf y}({\bf x},t)$ we introduce the averaged motions
${\bf y}^n = \phi_n \star {\bf y}$, where $\phi$ is a mollifier, and define ${\bf y}$ to be a singular limiting induced from
continuum solution (in short a {\it slic}-solution) if for every choice of the mollifier the smooth approximating family 
in the limit of small-scale averaging gives
\begin{equation}
\label{introapproxelas3d}
{\bf y}^n_{tt} - \div {\bf S} (\nabla {\bf y}^n ) =: {\bf f^n}  \to 0  \quad \mbox{in $\cD'$}.
\end{equation}
In this notion the precise form of regularizing mechanisms is not taken into account, instead it is enforced that they
act in a stable way amounting to averaging of the tested solution. Moreover, the energy equation for
\eqref{introapproxelas3d}, 
\begin{equation}
\del_t \left ( \frac{1}{2}  |{\bf y}_t^n |^2 + W( \nabla {\bf y}^n ) \right ) - \div \big ( {\bf y}_t^n \cdot {\bf S} (\nabla {\bf y}^n ) \big ) =
{\bf y}_t^n \cdot {\bf f^n} \, , 
\end{equation}
suggests that,
even though ${\bf f^n}  \to 0$ in $\cD'$,  the power of the "microscopic forces" ${\bf y}_t^n \cdot {\bf f^n}$ may well
have a non-trivial contribution in the limit, which needs to be calculated.

The definition of slic-solutions exploits the fact that  the momentum equation \eqref{eq:introelas3d}
is a second-order evolution.
It has certain conceptual analogies to the approach of weak asymptotic solution employed by
Danilov and Shelkovich \cite{DS05} in order to define solutions for special systems of conservation laws
involving delta shocks as asymptotic limits of smooth 
approximate solutions. There are also certain differences with the concept of weak asymptotic solution:
Slic-solutions are generated via 
averagings  of a candidate  discontinuous solution, and are natural in 
a context of discontinuous solutions  for second order evolution equations. The latter property
provides a mechanical intuition for the definition and suggests its name.

We examine the ramifications of this definition in three examples. In section \ref{sec-fracture} we consider the equations
of one-dimensional elasticity
\begin{equation}
\label{intro-elas1d}
y_{tt} - (\tau (y_x) )_x = 0
\end{equation}
and test a specific example of a motion,
\begin{equation}
\label{intro-ytx}
y(x,t) =
 \begin{cases}
\lambda x  \charf_{x < -\sigma t}
+ ( - t Y(0) + \alpha x)  \charf_{- \sigma t < x < 0} 
+ ( t Y(0) + \alpha x)  \charf_{0 < x < \sigma t}
            + \lambda x  \charf_{\sigma t < x} 
  & t > 0
  \\
  \lambda x & t < 0
  \end{cases}
\end{equation}
towards being a slic-solution. The example \eqref{intro-ytx}
is a counterpart  (in the one-dimensional case) of the cavitating solutions (for $d\ge 3$) of \cite{ps88,ps98} and corresponds
to a crack forming out of a homogeneously deformed state at the location $x=0$ at time $t=0$, in conjunction with
two outgoing Lax-shocks propagating at $x = \pm \sigma t$. It is more singular than the dynamic cavitating 
solution in \cite{ps88} as 
$y_x$ has a delta-mass at the origin. A definition of the notion of slic-solution is given
in Definition \ref{defslicconvx}, conditions are then provided under which \eqref{intro-ytx} is a slic-solution 
for \eqref{intro-elas1d} in Proposition \ref{propo-slic1d} ,
and the energy balance for the approximate solution of the crack is computed in Proposition \ref{propo-eb1d}. It turns out
that \eqref{intro-ytx} can be interpreted as a slic-solution, at the same time there is an energetic cost for creating the crack 
that is projected in the limiting energy balance equation \eqref{energybal}.

In section \ref{sec:cavitation}, we consider the dynamically cavitating solution  \eqref{eq:introradsym}, \eqref{eq:intrordef} 
with $r(s)$ as constructed in \cite{ps88} (see Lemma \ref{lem:psest} stating its properties). 
The notion of slic-solution adapted to that example is provided in Definition \ref{def:slic3d} 
and a natural definition for the energy is given in Definition \ref{def:energy3d}. The analysis is more cumbersome
and is based on detailed estimations of the layers of the approximate solution, but the results 
parallel those of the one-dimensional example. Namely, conditions are given under which the cavitating solutions
provide a slic-solution in Theorem \ref{lem:slic3d}, and the energy balance of the approximate solutions is
computed in Proposition \ref{lem:energy3d2}. Regarding the issue of uniqueness of entropy weak solutions,
it turns out that if the solution is construed as an entropy weak solution, then there is non-uniqueness. By contrast, if
the solution is construed in the slic-sense, then there is a contribution to the total energy in the process of forming the cavity
which results to the energy after the cavity formation being larger than before the cavitation 
(Proposition \ref{lem:increase}).  This indicates that
the notion of entropy weak solutions is inadequate when dealing with discontinuous solutions and in particular it cannot
account for the work needed to create the cavity. A more discriminating concept of solution has to be employed on
strong singularities, and the slic-solution concept is such a possibility.

The last example is the Riemann problem for the p-system in Lagrangean coordinates \eqref{eq:psys} 
with vacuum initial data \eqref{eq:vacid}. 
It has been conjectured that delta-shocks are needed to resolve the Riemann
problem with vacuum  \cite[Sec 9.6]{dafermosbook}, A proper definition for such solutions is not available at the hyperbolic
level, nevertheless solutions with delta shocks are constructed by using
viscous wave fans in the zero-viscosity limit ({\it e.g.} \cite{KK89}, \cite{E00}, \cite[Sec 9.8]{dafermosbook}). 
In section \ref{sec-vacuum}, we provide a definition for self-similar slic-solutions (Definition \ref{def:selfsimslic})
which is capable to define Riemann solutions with
delta-shocks at the hyperbolic level without recourse to a construction method (Theorem \ref{lem:vac:sol}) and 
to calculate their  mechanical energy  (Proposition \ref{lem:vac:energy}).


\section{A special motion exhibiting fracture in one-dimensional elasticity}
\label{sec-fracture}

The equation 
\begin{equation}
\label{elas}
y_{tt} = \tau ( y_x )_x \, , \quad x \in \R, \; t > 0
\end{equation}
describes longitudinal or shearing  motions $y(x,t)$ of one-dimensional elastic bars. We consider 
loading situations that the bar is subjected to homogeneous deformations far away, 
that is for some $r$ sufficiently large
$$
y (x,t) = \lambda \, x \, , \quad \mbox{ for $|x| > rt$}.
$$
Introducing the velocity $v = y_t$ and the strain $u = y_x$, the equation \eqref{elas} is expressed as a system of
conservation laws
\begin{equation}
\label{elasys}
\begin{aligned}
u_t &= v_x 
\\
v_t &= \tau(u)_x \, .
\end{aligned}
\end{equation}
In what follows we will assume that the stress function $\tau (u)$ satisfies the hypotheses
\begin{equation}
\label{hypo1}
\tau'(u) > 0 \, , \quad \tau''(u) < 0
\tag{$a_1$}
\end{equation}
\begin{equation}
\label{hypo2}
\tau(u) \to - \infty \quad \mbox{as $u \to 0$} \quad \mbox{and} \quad \int_1^u \tau(s) ds \to + \infty \quad \mbox{as $u \to 0.$}
\tag{$a_2$}
\end{equation}
Under \eqref{hypo1} the wave speeds $\lambda_{1,2} (u) = \pm \sqrt{ \tau'(u) }$ are real and \eqref{elas} is hyperbolic.
The hypothesis $\tau''(u) < 0$ is appropriate for an elastic material exhibiting softening elastic response
and plays an important role in the forthcoming analysis. 
The hypothesis \eqref{hypo2} is applicable in the case of longitudinal motions and 
is placed to exclude that a finite 
volume is compressed down to zero. In the sequel we will consider only tensile deformations
and this hypothesis will not play any significant role.
 (In the case of shearing motions $\tau (u)$ is defined for $u \in \R$ and \eqref{hypo2} is removed).

Smooth solutions of \eqref{elas} satisfy the additional conservation law
\begin{equation}
\label{energyeq}
\del_t \big ( \frac{1}{2} y_t^2 + W(y_x) \big )  - \del_x \big ( \tau(y_x) y_t \big ) = 0
\end{equation}
where 
$W(u) := \int_1^u \tau(s) ds$ is the elastic stored energy function. 
This equation captures the conservation of mechanical energy, 
$$
\frac{d}{dt} \int_a^b \frac{1}{2} y_t^2 + W(y_x) \, dx = \tau(y_x) y_t \Big |_{x=a}^b 
$$
namely, for each subdomain $(a,b)$, the rate of mechanical energy is balanced by the fluxes through the boundaries.
For elastic bars subjected to homogeneous deformations for $|x| > r$, the
mechanical energy  in $(-r, r)$ is conserved.

\subsection{Outline of the example}
We are interested in self-similar solutions $y = t Y \big ( \frac{x}{t} \big )$, with $Y = Y(\xi)$ a function of the
similarity variable $\xi = \frac{x}{t}$. The function $Y$ solves the equation
\begin{equation}
\label{selfelas}
\xi^2 Y'' = \tau ( Y' )' \, .
\end{equation}
Alternatively, we may introduce in \eqref{elasys} the self-similar ansatz $u = u \big ( \frac{x}{t} \big )$,
$v = v \big ( \frac{x}{t} \big )$ and recast it in the form of the Riemann problem
\begin{equation}
\label{selfelasys}
\begin{aligned}
-\xi u'  &= v'
\\
-\xi v' &= \tau(u)'  \, ,
\end{aligned}
\end{equation}
where $u = Y'$ and $v = Y - \xi u$.

If $\lambda$ is sufficiently large we expect that the elastic bar will break. Once the bar breaks, coherence is lost and
the hypothesis of continuum implicit in the derivation of \eqref{elas} is no longer valid.
Nevertheless, there will be a transition region from a range of loading where the model is valid
to a range where it loses validity,  and in some situations one may arguably give some meaning to equation \eqref{elas} 
past the transition regime. It is
precisely this transition at the onset of fracture that we wish to explore.

The scope of this section is to test a class of self-similar solutions describing the onset of fracture in one-space dimension
that are  generated by the function  $Y$ :
\begin{equation}
\label{frsol}
Y (\xi) =
\begin{cases}
\; \; \,  \lambda \xi  &  \; \; \; \xi < -\sigma  \\
-Y(0) + \alpha \xi  & -\sigma < \xi < 0  \\
\; \; \,  Y(0) + \alpha \xi  & \; \; \;  0 < \xi < \sigma  \\
\; \; \,  \lambda \xi    & \; \; \; \sigma < \xi  
\end{cases}
\end{equation}
where $\alpha$, $\lambda$, $\sigma$ and $Y(0)$ are positive parameters that satisfy $\lambda > \alpha$, $Y(0) > 0$,
 and are connected through the equations
\begin{align}
\label{kinem}
Y(0) &= (\lambda - \alpha) \sigma
\\
\label{rh}
\sigma &= \sqrt{ \frac{ \tau(\lambda) - \tau (\alpha)}{\lambda - \alpha} }.
\end{align}
The function $Y$ is continuous at $\xi = \pm \sigma$ but discontinuous at the origin $\xi = 0$.
The associated distributions $u$ and $v$ are determined via:
\begin{align}
u &= Y' = \bar u (\xi) + 2 Y(0) \delta_{\xi = 0} \quad 
&&\mbox{ where} \quad
\bar u (\xi) = \begin{cases} \alpha & 0 < \xi < \sigma \\ \lambda & \sigma < \xi \end{cases}
\quad \; \mbox{and $\bar u(-\xi) = \bar u(\xi)$}
\\
v &= Y - \xi u = \bar v (\xi) \quad 
&&\mbox{ where} \quad
\bar v (\xi) = \begin{cases} Y(0)  & 0 < \xi < \sigma \\ 0 & \sigma < \xi \end{cases}
\; \; \mbox{and $\bar v(-\xi) = - \bar v(\xi)$}.
\end{align}

In the sequel we propose a notion of solution according to which the function \eqref{frsol} may be
interpreted as a solution of \eqref{elas}.
Several preliminary remarks are in order:

\begin{itemize}
\item[(i)] The solution \eqref{frsol} has several common features and was in fact inspired by the dynamic
cavitating solutions in three-dimensions proposed in the important work by Spector and Pericak-Spector \cite{ps88}. 
One important difference is that the cavitating solutions in \cite{ps88} do not involve a delta measure for the strain.
In the present one-dimensional situation the solution  \eqref{frsol} may be thought as describing fracture, and accordingly
we will call it a crack.

\item[(ii)] For $\lambda$ fixed, there is a one parameter family of functions \eqref{frsol}
that satisfy  \eqref{kinem}, \eqref{rh}.
\item[(iii)] Note that $\lim_{t \to 0+} t Y \big ( \frac{x}{t} \big ) = \lambda x$ and that 
$\lim_{t \to 0+}  v \big ( \frac{x}{t} \big ) = 0$. Accordingly, at time $t=0$ the bar is in a configuration
of homogeneous deformation with strain $\lambda$ and at rest.
\item[(iv)] For positive times the function $y = t Y \big ( \frac{x}{t} \big )$ has a jump at the
origin associated to a crack that is moving according to $y(0\pm, t) = \pm t Y(0)$. The parameter $Y(0)$ stands for the
speed of the crack.
\item[(v)] The singularities at $\xi = \pm \sigma$ are shocks. The Rankine-Hugoniot conditions at $\xi = \sigma$,
$$
\begin{aligned}
-\sigma ( \lambda - \alpha) &= - Y(0)
\\
- \sigma ( - Y(0) ) &= \tau (\lambda) - \tau (\alpha) \, ,
\end{aligned}
$$
are satisfied due to \eqref{kinem}, \eqref{rh}. This shock belongs to the second characteristic family. In view
of hypothesis \eqref{hypo1},  
the Lax shock-admissibility criterion
$$
\lambda_2 (u-) = \sqrt{ \tau' (\alpha)} > \sigma = \sqrt{ \frac{ \tau(\lambda) - \tau (\alpha)}{\lambda - \alpha} }
> \sqrt{ \tau' (\lambda)} = \lambda_2 (u+)
$$
will be satisfied provided that $\alpha < \lambda$. A similar analysis indicates that the shock at $\xi = - \sigma$
belongs to the first characteristic family and
satisfies the Rankine-Hugoniot conditions and the Lax shock-admissibility criterion.

\item[(vi)] The distribution $(u, v)$ contains the delta measure $\delta_{\xi = 0}$ and therefore 
on one hand it is more singular than the usual solution of the Riemann problem, on the other hand 
one needs to give a meaning to $\tau (u)$. A direct computation, using 
\eqref{kinem} and the identity $- \xi \del_\xi \delta_{\xi = 0} = \delta_{\xi = 0}$ in $\cD'$,
shows that
$$
\begin{aligned}
- \xi u'  &= \sigma (\alpha - \lambda) \delta_{\xi = - \sigma} - 2 Y(0) \xi \del_\xi \delta_{\xi = 0} 
-\sigma (\lambda - \alpha) \delta_{\xi = \sigma}
\\
&= - Y(0) \delta_{\xi = - \sigma} + 2 Y(0) \delta_{\xi = 0}  - Y(0) \delta_{\xi = \sigma}
\\
&= v'
\end{aligned}
$$
in the sense of distributions. Hence, \eqref{selfelasys}$_1$ is satisfied. 
It remains to give a meaning to \eqref{selfelasys}$_2$.
Indeed, this is a main task of this section, and will be done 
through the notion of {\it singular limiting induced from continuum solution} that we introduce in the sequel.

\end{itemize}

\begin{figure}
    \centering 
   (a)\quad {\includegraphics[height=5cm,width=7cm]{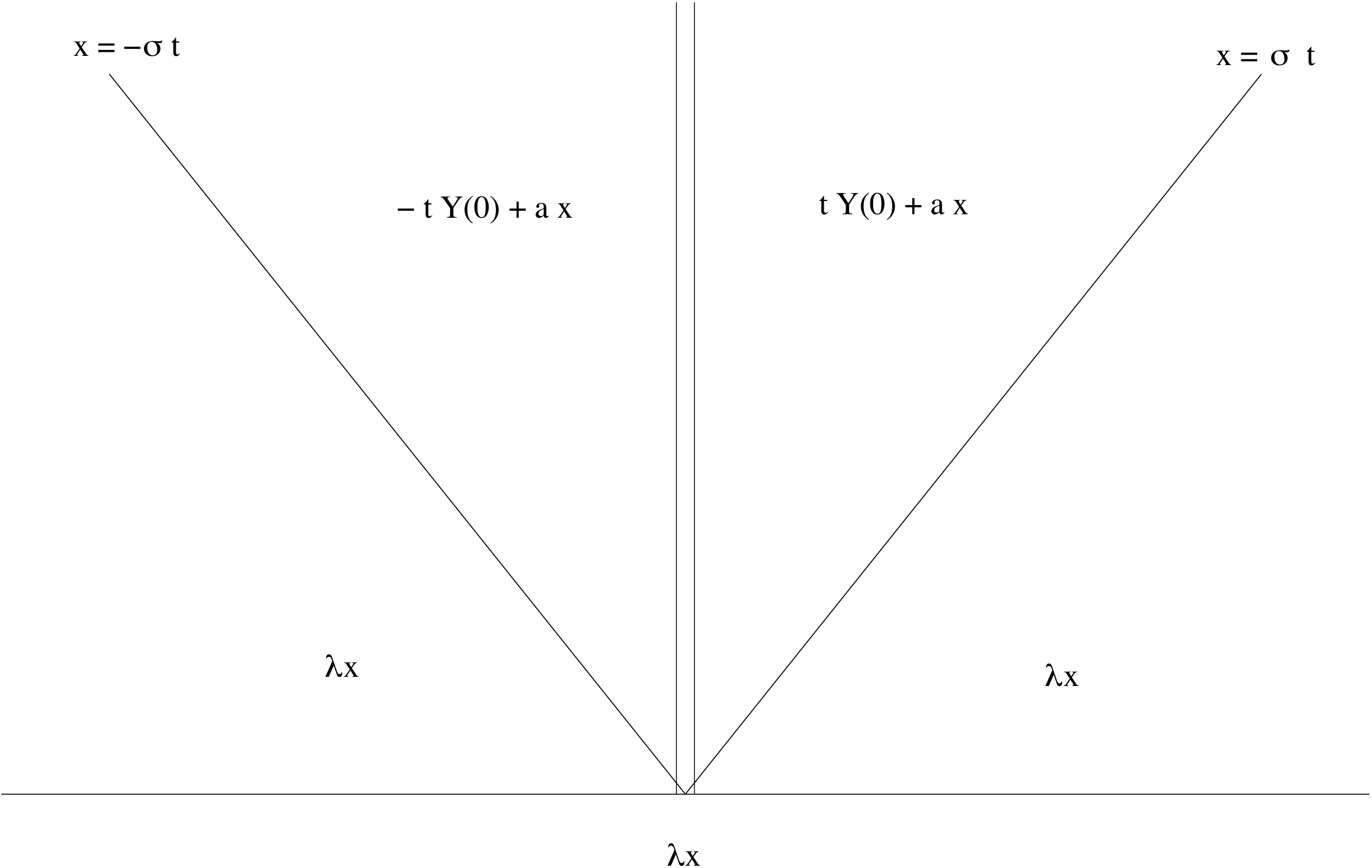}}
    \hspace{1.5cm}
   (b)\quad {\includegraphics[height=5cm,width=5cm]{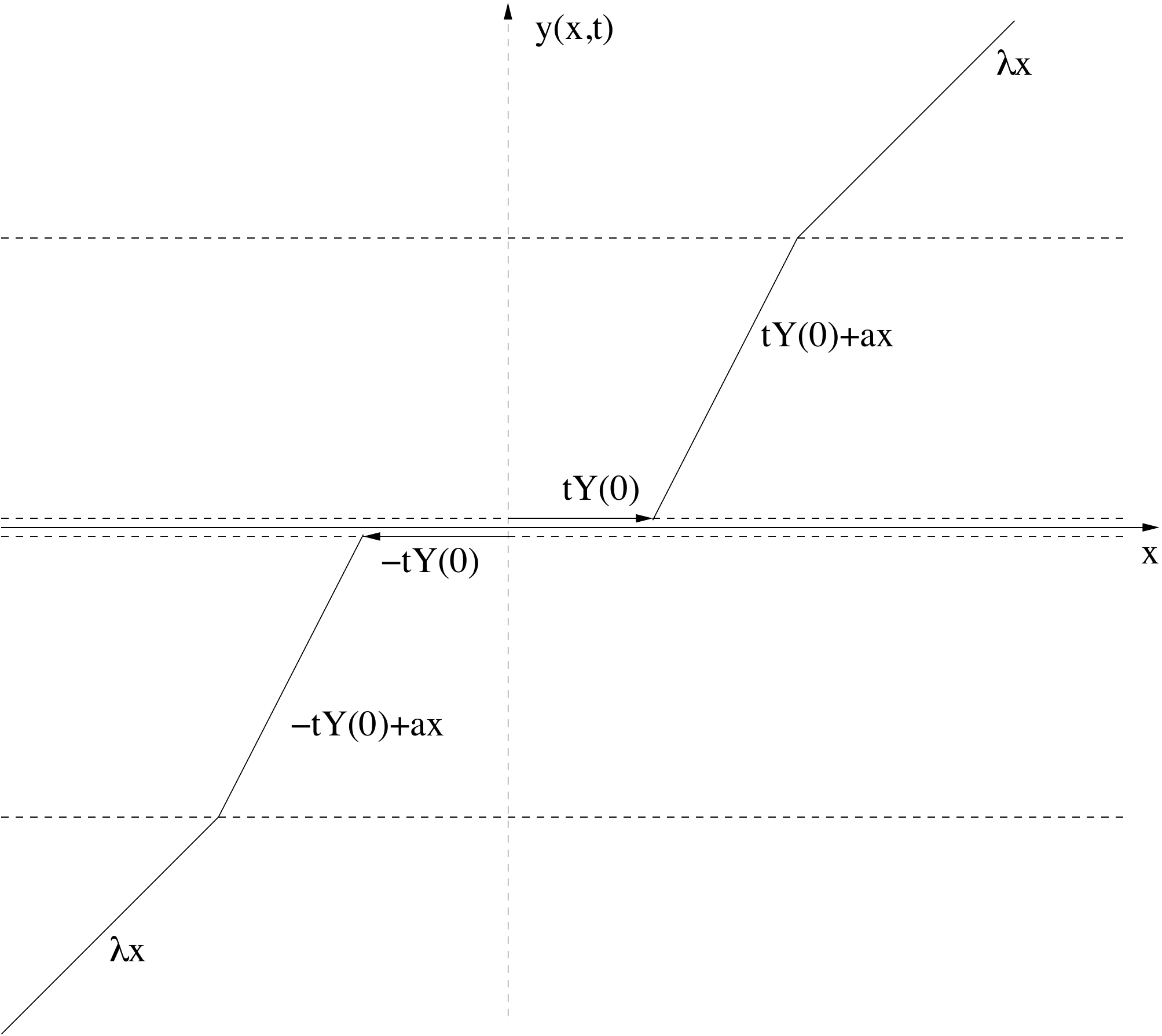}}    
  \caption{(a) cavitating solution in longitudinal motion ;  (b) shear band } 
  \label{fig1}
\end{figure}

In what follows we test $y = t Y \big ( \frac{x}{t} \big )$, with $Y = Y(\xi)$ defined in \eqref{frsol},
towards being a solution (interpreted in an appropriate sense) of \eqref{elas}.  It is expedient to extend the function $y$ for $t< 0$ by setting
$y = \lambda x$. 
The extended function still denoted by $y$ reads
\begin{equation}
\label{ytx}
y(x,t) =
 \begin{cases}
\lambda x  \charf_{x < -\sigma t}
+ ( - t Y(0) + \alpha x)  \charf_{- \sigma t < x < 0} 
+ ( t Y(0) + \alpha x)  \charf_{0 < x < \sigma t}
            + \lambda x  \charf_{\sigma t < x} 
  & t > 0
  \\
  \lambda x & t < 0
  \end{cases}.
\end{equation}
The parameters satisfy \eqref{kinem}, \eqref{rh} and $\lambda > \alpha$. Note that $y$ has a discontinuity
at $x=0$, but it is continuous at the shocks  $x = \pm \sigma t$ due to \eqref{kinem}. For longitudinal deformations
of an elastic bar, the motion \eqref{ytx} may be interpreted as a crack emerging at time $t=0$
from a homogeneously deformed state; the form of the solution in the $x-t$ plane is shown in Figure \ref{fig1}(a). The equation \eqref{elas} may also
be interpreted as describing elastic shear motions. In this context \eqref{ytx} represents a shear band emerging at $t=0$
and is sketched (in a $y(x,t) - x$ graph) in Figure \ref{fig1}(b).

For future reference we calculate certain properties for the solution: First, note that for $t > 0$
\begin{equation}
\label{prop1}
\begin{aligned}
y(x , t ) - y(x,0-) &= -t Y(0) \charf_{- \sigma t < x < 0}  + t Y(0) \charf_{0 < x < \sigma t} 
+ (\alpha - \lambda) x \charf_{ -  \sigma t  < x < \sigma t}
\\
| y(x , t ) - y( x, 0-) | &\le  \; t |Y(0)| + (\lambda - \alpha) \sigma t  \to 0  \qquad  \mbox{as $t  \to 0$}.
\end{aligned}
\end{equation}
Furthermore, in view of the kinematic compatibility assumption \eqref{kinem} and \eqref{prop1},
the distributional derivatives $\del_x y$ and $\del_t y$ are given by
\begin{equation}
\label{prop2}
\begin{aligned}
\del_x y &= 
\Big ( 2t Y(0) \delta_{x=0} + \lambda \charf_{x<-\sigma t}
+ \alpha \charf_{-\sigma t < x < \sigma t} + \lambda \charf_{\sigma t < x} \Big ) \charf_{t > 0} + \lambda \charf_{t < 0}
\\
\del_t y &=  \big ( -Y(0) \charf_{-\sigma t < x < 0} + Y(0) \charf_{0 < x < \sigma t} \big ) \charf_{t > 0}.
\end{aligned}
\end{equation}
Observe that $\del_t y (x,t) \to 0$ as $t \to 0$ for $x \ne 0$.

\subsection{Slic-solutions - definition}

There is a class of problems in material science where structures with discontinuous displacement fields
emerge via a material instability mechanism. Typical examples are development of cracks in fracture, 
cavitation in elastic response, or formation of shear bands in plastic deformations. Such problems lie at
the limits of applicability of continuum modeling and the usual concept of weak solutions is in
any case inadequate to describe these motions. Nevertheless, as the material transitions
from a regime where continuum modeling is applicable to a regime that the model has to be modified (or perhaps
atomistic modeling has to be employed), it is expected that at the interface both types of modeling have to apply
 in an intermediate regime. It is further expected that such structures should appear in a small parameter limit of more complex models
that incorporate "higher-order physics", and that their appearance occurs in a stable way.

\smallskip
The concept of  {\it singular limiting induced from continuum solution} (in short {\it slic}-solution) is an attempt
to give meaning to such discontinuous solutions. It is presented  here at the level of the equations of
nonlinear elasticity. Roughly speaking a discontinuous solution of \eqref{elas} will be a slic-solution
if it can be obtained as the limit of approximate smooth solutions that are an averaging of $y$.
More precisely:

\begin{definition}
\label{defslic}
Let ${\bf y} \in L^1_{loc} ( \cQ)$ where $\cQ$ is an open domain in space-time. 
Given a mollifier $\varphi \in C_c^\infty (\R^d \times \R )$, with $\varphi \ge 0$, $\supp \varphi \subset \{ |{\bf x}| < 1, t \in (-1,1) \} $, 
$\iint \varphi d{\bf x} dt= 1$, we set  $\varphi_n =n^{d+1}  \varphi  ( n {\bf x} ,  n t  )$.
For any $\cO \Subset \cQ$, $\cO$ compactly embedded in $\cQ$, the averaged function
\begin{equation}
\label{averaging}
{\bf y}^n ({\bf x}, t) =  \varphi_n  \underset{{\bf x},t}{\star}  {\bf y} = \iint \varphi_n ( {\bf x} - {\bf z}, t - \tau ) 
{\bf y}( {\bf z},\tau) \, d{\bf z} \, d\tau
\end{equation}
is well defined for $({\bf x},t) \in \cO$ and $\frac{1}{n} < \dist ( \cO , \partial \cQ)$.
We will say that ${\bf y}$ is a {\it singular limiting induced from continuum solution} of \eqref{eq:introelas3d} if,
for any mollifier $\varphi,$ any $\cO \Subset \cQ$ and for $\psi \in C_c^\infty ( \cO \, ; \,  R^d )$,
\begin{equation}
\label{slic}
\iint  {\bf y}^n \psi_{tt} + {\bf S} (\nabla {\bf y}^n) : \nabla \psi  \, d {\bf x} dt  \to 0 \quad \mbox{ as $n \to \infty$.}
\end{equation}
\end{definition}

Several remarks are in order. This definition incorporates as a requirement
that the singular solution results by a smooth averaging process, but it does not contain any additional
information on the approximation process that might induce a mechanism for selecting among singular solutions.
One can easily  check that the usual concept of weak solutions (containing shocks) 
of the equations of elasticity are also slic-solutions.  Finally, the concept of slic-solution is quite natural
when continuum modeling is viewed as resulting from the discrete, especially in situations where the distance
between lattice elements is large but nevertheless matter is present. In that sense it is more natural 
in a context of shear bands than in a context of cavitation after a hole has fully formed.

The definition of {\it slic}-solution is  cumbersome to use in practice. Nevertheless, it
provides a theoretical framework to define such singular solutions and to calculate the ramifications
on the energy. A difficulty occurs because in the general setting of Definition \ref{defslic}
averaging in both space and time is required . For hyperbolic problems we expect that
averaging in space will  also induce averaging in time and in many settings the double averaging might be avoided.
Also, in practice we expect to check the definition at only local singularities what simplifies calculations.
In particular when checking it for self-similar solutions averaging in space only is adequate.

For certain special systems of conservation laws, like the system of pressureless gas dynamics and the so called
Keyfitz-Kranzer system, there has been suggested to use  solutions that take values in measures,
a literature going by the name of delta-shocks ({\it e.g.} \cite{KK89} \cite{E00}, \cite[Sec 9.6]{dafermosbook} and references therein). 
Danilov and Shelkovich \cite{DS05} have proposed to define delta shocks as asymptotic limits of smooth 
approximate solutions, a process that they call weak asymptotic solution. This definition has obvious similarities to
the concept of  slic-solution presented here. There are also two differences: first, the notion of slic-solution
 is connected via averaging to a possibly discontinuous "solution candidate", second, 
slic-solutions emerge  in a context of discontinuous solutions  for second order evolution equations, 
what provides a mechanical intuition for the definition and suggests the name.

\subsection{The example of 1-d fracture as a slic-solution}

Next, we turn to the example of the fracturing solution \eqref{ytx} with the objective to test it as a slic-solution 
for the equations \eqref{elas}.
As already mentioned the goal is to give a meaning to the solution at $x=0$. To simplify calculations 
we will introduce a variant
of Definition \ref{defslic}  that uses a mollification in the space variable $x$ only:

\begin{definition} 
\label{defslicconvx}
Let $y \in L^\infty_{loc} \big ( (-\infty ,\infty) \, ; \, L^1_{loc} (\R) \big )$
 satisfy for some $\eps  > 0$
the monotonicity condition:
\begin{equation}
\label{hypomc}
y(x_1,t)-y(x_2,t)>  \eps (x_1-x_2)  \quad \mbox{for $x_1,x_2,t \in \R$ with $x_1>x_2$.}
\tag{$mc$}
\end{equation}
For  $\phi$ a mollifier,
 $\phi \in C^\infty_c (\R)$, $\phi \ge 0$,  $\supp \phi \subset  B_1$ (the ball of radius 1),
$\int \phi = 1$, we let $\phi_n = n  \phi  ( n x   )$ and  define the averaged function
\begin{equation}
\label{convx}
y^n (x,t) = \phi_n \underset{x}{\star}  y = \int \phi_n (x-z) y(z,t) dz.
\end{equation}
The function $y$ is called a singular limiting induced from continuum (slic)  solution of \eqref{elas} provided
for any symmetric mollifier
\begin{equation}
\label{eqdefslic}
 \int_\R \int_\R y^n \psi_{tt} + \tau(  y^n_x ) \psi_x \, dx dt \rightarrow 0
 \end{equation}
as $n \to \infty$ for $ \psi \in C_c^\infty (\R \times \R).$
\end{definition}


To justify that the left hand side of \eqref{eqdefslic} is well defined, observe that in view of \eqref{convx}
$$
\begin{aligned}
|\del_x y^n (x,t)| &\le n^2 \| \phi' \|_{C^0}  \int_{|z-x| < \frac{1}{n}}  |y(z,t)| dz \, ,
\end{aligned}
$$
and thus for $y$ of class $L^\infty_{loc} \big ( (-\infty ,\infty) \, ; \, L^1_{loc} (\R) \big )$, the derivative
$\del_x y^n \in L^\infty_{loc} ( (-\infty ,\infty) \times \R)$. The monotonicity condition placed on ${\bf y}$ ensures
that $\del_x {\bf y}^n  \ge \eps $ for all $n.$ 
The definition in the form stated applies to longitudinal motions.  
In the case that \eqref{elas} models shearing motions then hypothesis \eqref{hypomc}
is  removed and the condition $\del_x {\bf y}^n \ge \eps > 0$ is no longer necessary.
Definition \ref{defslic} using the convolution
in $x$ formula \eqref{convx}  relies on the special form of \eqref{elas} and is easier to use in practice. 

We turn now to the example \eqref{ytx} with $\alpha$, $\lambda$, $Y(0)$ satisfying \eqref{kinem}, \eqref{rh}
and do some preparatory computations.
For $t > 0$,
$$
\begin{aligned}
u^n &:= \del_x y^n  =  \int_\R \del_x \phi_n (x-z) \, y(z,t) dz
\\
&= - \int_\R \del_z  ( \phi_n (x-z)  ) \Big [ \lambda z \charf_{z < -\sigma t} +
(-t Y(0) + \alpha z) \charf_{-\sigma t < z < 0} + (t Y(0) + \alpha z) \charf_{0 < z < \sigma t} + 
\lambda z \charf_{\sigma t < z} \Big ] dz
\\
&=
t \big ( \sigma (\lambda - \alpha) - Y(0) \big ) \phi_n (x + \sigma t) + 2 \phi_n (x) t Y(0) 
+ t \big ( \sigma (\lambda - \alpha) - Y(0) \big ) \phi_n (x - \sigma t)
\\
&\quad 
+ \lambda \int_{-\infty}^{-\sigma t} \phi_n (x-z) dz 
+ \alpha \int_{-\sigma t}^{\sigma t} \phi_n (x-z) dz
+ \lambda \int_{\sigma t}^{\infty} \phi_n (x-z) dz
\\
&\stackrel{\eqref{kinem}}{=}
2 \phi_n (x) t Y(0) + \lambda \int_{-\infty}^{-\sigma t} \phi_n (x-z) dz 
+ \alpha \int_{-\sigma t}^{\sigma t} \phi_n (x-z) dz
+ \lambda \int_{\sigma t}^{\infty} \phi_n (x-z) dz
\end{aligned}
$$
and 
$$
\begin{aligned}
v^n := \del_t y^n &= \int \phi_n (x-z) \del_t y(z,t) dz
\\
&= - Y(0) \int_{-\sigma t}^0 \phi_n (x-z) dz + Y(0) \int_0^{\sigma t} \phi_n (x-z) dz.
\end{aligned}
$$
Since $\alpha < \lambda$ and $Y(0) > 0$, we have the bounds
\begin{equation}
\label{uvbound}
\begin{aligned}
\alpha \le u^n  &\le \lambda + 2 t Y(0) \phi_n (x)
\\
|v^n| &\le Y(0)
\end{aligned}
\end{equation}
In addition, we compute
\begin{align}
\label{yx}
\del_x y^n &= 
\begin{cases}
2 \phi_n (x) t Y(0) + \lambda \int_{-\infty}^{-\sigma t} \phi_n (x-z) dz 
+ \alpha \int_{-\sigma t}^{\sigma t} \phi_n (x-z) dz
+ \lambda \int_{\sigma t}^{\infty} \phi_n (x-z) dz
&  t > 0  \\
\lambda & t < 0  \\
\end{cases}
\\
\label{yvn}
\del_t y^n &=
\begin{cases}
- Y(0) \int_{-\sigma t}^0 \phi_n (x-z) dz + Y(0) \int_0^{\sigma t} \phi_n (x-z) dz
&  t > 0  \\
0 & t < 0  \\
\end{cases}
\\
\label{yan}
\del_{t}^2 y^n &=
\begin{cases}
- Y(0)  \sigma \phi_n (x + \sigma t ) + Y(0) \sigma  \phi_n (x- \sigma t)
&  t > 0  \\
0 & t < 0  \\
\end{cases}.
\end{align}
It is a tedious but straightforward calculation to check that, mainly as a consequence of \eqref{kinem}, 
the regularity of $y^n(x,t)$ is of class 
 $ C^2 (\R^2)$ (but $\del^3_t y$ has a jump
discontinuity at $t=0$). One also easily checks that
\begin{align*}
\del_x y^n  &\stackrel{\cD'_{x,t}}{\longrightarrow}  
\Big ( 2t Y(0) \delta_{x=0} + \lambda \charf_{x<-\sigma t}
+ \alpha \charf_{-\sigma t < x < \sigma t} + \lambda \charf_{\sigma t < x} \Big ) \charf_{t > 0} + \lambda \charf_{t < 0}
\\
\del_t y^n  &\stackrel{\cD'_{x,t}}{\longrightarrow}    \big ( -Y(0) \charf_{-\sigma t < x < 0} + Y(0) \charf_{0 < x < \sigma t} \big ) \charf_{t > 0}
\end{align*}
which should be compared to \eqref{prop2}.

\medskip

Next, we  test whether $y$ in \eqref{ytx} provides a slic-solution for \eqref{elas}.

\begin{proposition}
\label{propo-slic1d}
If \eqref{kinem}, \eqref{rh} are satisfied and 
\begin{equation}
\label{hypo3}
L : = \lim_{ u \to \infty} \frac{\tau(u) }{u} = 0
\tag{$a_3$}
\end{equation}
then $y$ defined in \eqref{ytx} is a slic-solution for \eqref{elas}.
\end{proposition}

\proof 
Due to  \eqref{hypo1}, the limit \eqref{hypo3} is always well defined:
Indeed, since $\tau(u)$ is monotone increasing, either $\tau_\infty := \lim_{u\to \infty} \tau(u) < \infty$ and $L= 0$,
 or $\tau_\infty = \infty$ and $L = \lim_{u \to \infty} \tau' (u)$  exists by  \eqref{hypo1}.
 
For $y$ to be a slic-solution,  its regularization $y^n = \phi_n \underset{ x}\star y$ should satisfy
the approximate equation
\begin{equation}
\label{appeq}
\del^2_t y^n - \del_x \tau (\del_x y^n ) =: f^n
\end{equation}
with $f^n \to 0$ in $\cD'$. As $y^n \in C^2 (\R^2)$ the equation \eqref{appeq} is satisfied  in a classical sense,
and it will be viewed as defining $f^n$. 
Let $\psi \in C_c^\infty ( \R \times \R)$ be a test function, and consider the distribution $< f^n , \psi >$
associated to $f^n$. The time-domain is split  in two parts:  $0< t < \frac{2}{n\sigma}$ and $t > \frac{2}{n\sigma}$. The
 threshold  $\frac{2}{n\sigma}$ is selected so that the waves 
in the approximate solution $y^n$ are separated and no longer interact for $t > \frac{2}{n\sigma}$. 
\begin{equation}
\label{fn}
\begin{aligned}
< f^n , \psi > &= \int_{-\infty}^\infty \int_{\R} (\del_t^2 y^n ) \psi + \tau(\del_x y^n ) \del_x \psi \, dx dt
\\
&= \int_{\frac{2}{n\sigma}}^\infty \int_\R \del_t^2 y^n  \psi dx dt 
+ \int_{\frac{2}{n\sigma}}^\infty \int_{\R \backslash (-\frac{1}{n}, \frac{1}{n})}  \tau(\del_x y^n) \del_x \psi \, dx dt
\\
&\quad+ \int_{\frac{2}{n\sigma}}^\infty \int_{(-\frac{1}{n}, \frac{1}{n})} \tau(\del_x y^n) \del_x \psi \, dx dt
+ \int_0^{\frac{2}{n\sigma}} \int_{\R} \del_t^2 y^n  \psi + \tau(\del_x y^n) \del_x \psi \, dx dt
\\
&= : J^n + I_1^n + I_2^n + E^n \, .
\end{aligned}
\end{equation}

We proceed to estimate the right hand side of \eqref{fn} starting with the term $E^n$. Observe that
\begin{equation}
\label{bounduinit}
\alpha \le \del_x y^n \le \lambda + \frac{4}{\sigma} Y(0) \| \phi \|_0 \, , \quad \mbox{for $0 < t < \frac{2}{n\sigma}$, $x \in \R$}
\end{equation}
and thus
\begin{equation}
\label{En}
\begin{aligned}
|E^n| &\le 
\int_0^{\frac{2}{n\sigma}} \int_{\R} 
Y(0) \sigma (  \phi_n (x+\sigma t) + \phi_n (x -\sigma t) ) |\psi| + |\tau(\del_x y^n) | |\del_x \psi |\, dx dt
\\
&\le Y(0)  \frac{4}{n} \| \psi \|_0  \; + \;
\max\big\{ | \tau(\alpha)|, \big | \tau \big ( \lambda + \frac{4}{\sigma} Y(0) \| \phi \|_0 \big ) \big | \big\} \int_0^{\frac{2}{n\sigma}} \int_{\R}  |\del_x \psi |\, dx dt
\\
&\to 0  \quad \mbox{as $n \to \infty$}.
\end{aligned}
\end{equation}
Using \eqref{yan}, we next obtain
\begin{equation}
\label{jn}
\begin{aligned}
J^n &= \int_{\frac{2}{n\sigma}}^\infty \int_\R \del_t^2 y^n  \psi dx dt 
\\
&= \int_{\frac{2}{n\sigma}}^\infty \int_\R \big ( -Y(0) \sigma \phi_n (x+\sigma t) + Y(0) \sigma \phi_n (x -\sigma t) \big ) \psi(x,t) \, dx dt
\\
&\to  \int_0^\infty -Y(0) \sigma \psi (-\sigma t, t) + Y(0) \sigma \psi(\sigma t, t) \, dt.
\end{aligned}
\end{equation}
Using \eqref{yx}, and the bound $\alpha \le \del_x y^n \le \lambda$ for $x \in \R \backslash (-\frac{1}{n}, \frac{1}{n})$, $t > \frac{2}{n\sigma}$,
we have
\begin{equation}
\label{i1n}
\begin{aligned}
I_1^n &= \int_{\frac{2}{n\sigma}}^\infty \int_{\R \backslash (-\frac{1}{n}, \frac{1}{n})}  \tau(\del_x y^n) \del_x \psi \, dx dt
\\
&\to \int_0^\infty \Big ( \int_{-\infty}^{-\sigma t} \tau (\lambda) \psi_x dx + \int_{-\sigma t}^0 \tau (\alpha) \psi_x dx
+ \int_0^{\sigma t} \tau(\alpha) \psi_x dx + \int_{\sigma t}^\infty \tau(\lambda) \psi_x dx \Big ) dt.
\end{aligned}
\end{equation}

Finally, we compute the contribution of the origin:
\begin{equation}
\label{i2n}
\begin{aligned}
I_2^n &= \int_{\frac{2}{n\sigma}}^\infty \int_{-\frac{1}{n}}^{\frac{1}{n}} \tau \big ( \alpha + 2 \phi_n (x) t Y(0) \big ) \psi_x (x,t) \, dx dt 
\\
&= \int_{\frac{2}{n\sigma}}^\infty \int_{-\frac{1}{n}}^{\frac{1}{n}} 
\frac{ \tau \big ( \alpha + 2 \phi_n (x) t Y(0) \big ) }{\alpha + 2 \phi_n (x) t Y(0)}
\big ( \alpha + 2 \phi_n (x) t Y(0) \big ) \psi_x (x,t) \, dx dt 
\\
&= \int_{\frac{2}{n\sigma}}^\infty \int_{-1}^{1}  \frac{ \tau \big ( \alpha + 2 n \phi(z) t Y(0) \big ) }{ \alpha + 2 n \phi(z) t Y(0) }
\frac{ \alpha + 2 n \phi(z) t Y(0) }{n} \psi_x ( \frac{z}{n}, t) \, dz dt
\\
&\to 2 L Y(0) \int_0^\infty \psi_x (0,t) t \, dt.
\end{aligned}
\end{equation}

Combining \eqref{fn} with \eqref{En}, \eqref{jn}, \eqref{i1n}, \eqref{i2n} and \eqref{rh}, we see that
\begin{equation}
\label{slic1d}
\begin{aligned}
<f^n , \psi> &= J^n + I_1^n + I_2^n  + E^n
\\
&\to 
\int_0^\infty \Big (  [ \tau(\lambda) - \tau (\alpha) - \sigma Y(0) ] \psi (-\sigma t, t) 
- [ \tau(\lambda) - \tau (\alpha) - \sigma Y(0) ] \psi (\sigma t, t) 
\\
&\qquad \qquad \qquad + 2 Y(0) t L \del_x \psi(0,t) \Big )  dt
\\
&=  - 2 Y(0) L <  \del_x \delta_{x=0} , t \psi >
\end{aligned}
\end{equation}
and the latter is zero if and only if $L=0$. \qed

%
%
%

\subsection{The energetic cost for the creation of the crack}

The function $y^n$ is smooth and satisfies the approximate equation \eqref{appeq} and the
associated energy identity 
$$
\del_t \left ( \frac{1}{2} (v^n)^2 + W( y^n_x ) \right )  - \del_x ( v^n \tau ( y^n_x  ) ) = f^n v^n.
$$
Even though $f_n  \stackrel{\cD'_{x,t}}{\longrightarrow}   0$ the
product $f^n v^v$ may yield in the limit a nontrivial contribution to the energy. In this subsection we 
compute that contribution.

To consider the energy balance, fix a domain  $B = (-r, r)$ containing the entire wave fan of the approximate solution
\eqref{convx} at time $t$. For instance, if we are interested in times $0 < t \le T$, we select $r > \sigma T + 1$.
The velocity at the boundary of such a domain vanishes,
$v^n \big |_{\del B} = 0$. The total energy of the wave fan
\begin{equation}
\label{toten}
E [ y^n ; B  ] = \int_B \frac{1}{2} (v^n)^2 + W( \del_x y^n)  \, dx
\end{equation}
evolves according to the energy balance equation
$$
\frac{d}{dt} \int_B \frac{1}{2}(v^n)^2 + W( \del_x y^n)  dx = \int_B f^n v^n \,dx \, .
$$

\begin{lemma} 
If $v^n \big |_{\del B} = 0$, then  for $t>0$ the (rate of) change of total energy $E[y^n ; B ]$ is calculated via
\begin{equation}
\label{energy}
\begin{aligned}
\frac{d}{dt} E [ y^n ; B ] &=
\Big ( Y(0)^2 \sigma + 2 \sigma [W(\alpha) - W(\lambda)]   + 2 Y(0) \tau (\alpha) \Big ) \charf_{t > \frac{2}{n \sigma}}
\\
&\quad + \charf_{t > \frac{2}{n \sigma}} \left ( \int_{-\frac{1}{n}}^{\frac{1}{n}} \tau(\alpha + 2 \phi_n (x) t  Y(0)) 2Y(0) \phi_n (x) dx -2 Y(0) \tau(\alpha) 
\right )
\\
&\quad + e^n \charf_{0 < t < \frac{2}{n \sigma}} 
\end{aligned}
\end{equation}
where $e^n$ uniformly bounded independent of $n$.
\end{lemma}

\proof
We split the time-domain in two parts: the region $t > \frac{2}{n\sigma}$ where the energy dissipation
will be computed exactly, and the region $0< t < \frac{2}{n\sigma}$ where
it will be estimated. The choice of the threshold is so that for  $t > \frac{2}{n\sigma}$ the waves
in the approximate solution $y^n$ are separated and no longer interact.

First we consider the range $t > \frac{2}{n\sigma}$.
Using the notation $v^n = \del_t y^n$, $u^n = \del_x y^n$ and 
\eqref{yx}, \eqref{yvn}, \eqref{yan}, we compute
$$
\begin{aligned}
\del_t \int_B &\frac{1}{2} (v^n)^2 + W( \del_x y^n)  \, dx =
\int_B v^n \del_t v^n + \tau(u^n) \del_t u^n dx
\\
&= \int_{- \sigma t - \frac{1}{n}}^{- \sigma t  + \frac{1}{n}}  \big ( -Y(0) \sigma v^n + \sigma (\alpha - \lambda) \tau( u^n) \big ) \phi_n (x+\sigma t) \, dx
+ Y(0) \tau(\alpha)
\\
&\quad 
+  \int_{ - \frac{1}{n}}^{ \frac{1}{n}}    \tau (u^n) 2 Y(0) \phi_n (x) \, dx
- 2 Y(0) \tau(\alpha)
\\
&\quad 
+  \int_{ \sigma t - \frac{1}{n}}^{\sigma t  + \frac{1}{n}}  \big ( Y(0) \sigma v^n + \sigma (\alpha - \lambda) \tau( u^n) \big ) \phi_n (x - \sigma t) \, dx
+ Y(0) \tau(\alpha)
\\
&= :  \mu^n_{-\sigma}  + p^n_c + \mu^n_{\sigma}
\end{aligned}
$$

Using \eqref{yx} and \eqref{yvn} we obtain
$$
\begin{aligned}
\mu^n_{\sigma} 
&=  Y(0) \tau(\alpha) +   Y(0)^2 \sigma \int_{\sigma t - \frac{1}{n}}^{\sigma t 
+ \frac{1}{n}}  \int_0^{\sigma t} \phi_n(x-z) \phi_n (x-\sigma t) dz dx
\\
&\quad +  \int_{\sigma t - \frac{1}{n}}^{\sigma t + \frac{1}{n}}
\tau \Big ( \alpha \int_0^{\sigma t} \phi_n (x-z) dz + \lambda \int_{\sigma t}^\infty \phi_n (x-z) dz \Big )
(\alpha - \lambda) \sigma \phi_n (x-\sigma t) dx
\\
&=: Y(0) \tau(\alpha) + Y(0)^2 \sigma A^n  +    B^n.
\end{aligned}
$$

Observe now that
$$
\begin{aligned}
A^n &=  \int_{\sigma t - \frac{1}{n}}^{\sigma t 
+ \frac{1}{n}} \int_{-\infty}^{\sigma t} \phi_n(x-z) \phi_n (x-\sigma t) dz dx
\\
&
= \int_{- \frac{1}{n}}^{ \frac{1}{n}} \int_{-\infty}^{0} \phi_n(\bar x- \bar z) \phi_n (\bar x) d\bar z d\bar x
\\
&
= \int_{- \frac{1}{n}}^{ \frac{1}{n}} \int_{\bar x}^\infty \phi_n (u) \phi_n (\bar x)  du d\bar x 
\\
&= - \int_{- \frac{1}{n}}^{ \frac{1}{n}}  \frac{1}{2} \frac{d}{d\bar x} 
   \left ( \int_{\bar x}^\infty \phi_n (u) du \right )^2 d\bar x
= \frac{1}{2}
\end{aligned}
$$
and
$$
\begin{aligned}
B^n &= \int_{\sigma t - \frac{1}{n}}^{\sigma t + \frac{1}{n}}
\tau \Big ( \alpha + (\lambda - \alpha)  \int_{\sigma t}^\infty \phi_n (x-z) dz \Big )
[ - (\lambda - \alpha) \sigma] \phi_n (x-\sigma t) dx
\\
&= \int_{ - \frac{1}{n}}^{\frac{1}{n}}
\tau \Big ( \alpha + (\lambda - \alpha)  \int_{0}^\infty \phi_n (\bar x-\bar z) d\bar z \Big )
[ - (\lambda - \alpha) \sigma] \phi_n (\bar x) d\bar x
\\
&= (-\sigma) \int_{ - \frac{1}{n}}^{\frac{1}{n}} \frac{d}{d\bar x} W \Big ( \alpha + (\lambda - \alpha) 
\int_{-\infty}^{\bar x} \phi_n (u) du \Big ) d\bar x
\\
&= (-\sigma) \left [ W \Big ( \alpha + (\lambda - \alpha)  \int_{-\infty}^{\frac{1}{n}} \phi_n (u) du \Big )
- W \Big ( \alpha + (\lambda - \alpha)  \int_{-\infty}^{-\frac{1}{n}} \phi_n (u) du \Big ) \right ]
\\
&= \sigma ( W(\alpha ) - W(\lambda) ).
\end{aligned}
$$
Combining the above relations and performing a similar computation for the wave at the backward moving shock 
gives for the rate of dissipation of the two shocks:
\begin{alignat}{3}
\label{mu1}
\mu^n_\sigma  &= \mu_{\sigma} \; \;  = -\sigma \big [ -\frac{1}{2} Y(0)^2 + W(\lambda) - W(\alpha) \big ] + Y(0) \tau(\alpha)  
\quad &\mbox{$t >  \frac{2}{n\sigma}$},
\\
\label{mu2}
\mu^n_{-\sigma} &= \mu_{-\sigma} =  \; \;   \sigma \big [ \frac{1}{2} Y(0)^2 +W(\alpha)   - W(\lambda) \big ] + Y(0) \tau(\alpha) 
\quad &\mbox{$t >  \frac{2}{n\sigma}$}.
\end{alignat}
In the range $t > \frac{2}{n\sigma}$, the term $p^n_c$ is computed by
\begin{equation}
\label{enhole}
p^n_c =   \Big [ \int_{- \frac{1}{n}}^{\frac{1}{n}} \tau(\alpha + 2 \phi_n (x) t Y(0)) 2Y(0) \phi_n (x) dx
- 2 Y(0) \tau (\alpha) \Big ]
\end{equation}
and describes the (rate of)  work of the force $f^n$ at the crack thus determining the cost of energy for opening the 
hole.

Next, we focus on the range $0 < t < \frac{2}{n\sigma}$ and the term
\begin{equation}
\begin{aligned}
e^n (t)  &:= \int_B v^n \del_t v^n + \tau(u^n) \del_t u^n dx
\\
&= \int_B v^n  \Big ( -Y(0) \sigma \phi_n (x+\sigma t) + Y(0) \sigma \phi_n (x -\sigma t)   \Big )  
\\
&\quad  + \tau (u^n)  \Big ( 2 Y(0) \phi_n (x) +  \sigma (\alpha - \lambda) \big (  \phi_n (x+\sigma t) + \phi_n (x -\sigma t) \big ) \Big ) \, dx.
\end{aligned}
\end{equation}
Using  \eqref{uvbound} and \eqref{bounduinit}, we see that  $e^n \charf_{0<t < \frac{2}{n\sigma}}$ is uniformly bounded
independently of $n$.
\qed

The wave fan consists of three waves, two shocks located at $\xi = -\sigma$ and 
$\xi = \sigma$ and the delta-singularity associated to the crack and located at $\xi = 0$.
As already mentioned both shocks satisfy the Lax-admissibility criterion.  The standard Riemann problem 
theory, see \cite[Sec 8.5]{dafermosbook},  suggests that the dissipation at the shocks 
$\mu_\sigma$ and $\mu_{-\sigma}$ is computed by the formulas \eqref{mu1}, \eqref{mu2} and also that
$$
\mu_\sigma < 0 \, , \; \mu_{-\sigma} < 0 \, .
$$
The energy balance is thus expressed in the form  \eqref{energy} where $e^n$ stands for an error term,
$\mu_\sigma$ and $\mu_{-\sigma}$ is the dissipation at the two shocks and $p_c^n$ given by \eqref{enhole}
is the work of the force $f^n$ at the crack and describes the cost of energy for opening the 
hole. The limiting contribution of $p^n_c$ may be calculated 
$$
\begin{aligned}
p^n_c &=  \int_{-\frac{1}{n}}^{\frac{1}{n}} \tau(\alpha + 2 n \phi (n x) t Y(0)) 2Y(0) n \phi (n x) dx  - 2 Y(0) \tau (\alpha)
\\
&=  \int_{-1}^1 \tau(\alpha + 2 n \phi (z) t Y(0)) 2Y(0)  \phi (z) dz  - 2 Y(0) \tau (\alpha) 
\\
&\to
\begin{cases}
\infty &\mbox{if \quad  $\lim_{u \to \infty} \tau(u) = \infty$ } \\
2 ( \tau_\infty - \tau(\alpha) ) Y(0) &\mbox{if \quad $\lim_{u \to \infty} \tau(u) =: \tau_\infty < \infty$ }.\\
\end{cases}
\end{aligned}
$$
Accordingly, the cost of opening a cavity is infinite when  $\tau_\infty = \infty$. By contrast, if $\tau_\infty < \infty$
the energetic cost of opening the crack is finite and it is  conceivable that a fracture can appear.

We define the energy of the whole wave fan to be the limit of the energy of the approximate solution
\begin{equation}
E [ y ; B] = \lim_{n \to \infty} E [ y^n ; B ] \, .
\end{equation}
The last proposition shows that in the latter case the total energy after the crack is larger than the energy 
before the crack opens.

\begin{proposition}
\label{propo-eb1d}
If $\lim_{u \to \infty} \tau(u) =: \tau_\infty < \infty$ then the total energy of the approximate solution
 $E[y^n ; B]$ of the wave fan
satisfies
\begin{equation}
\label{energyfan}
\begin{aligned}
\frac{d}{dt} E [ y^n ; B] &= \left (  \mu_{-\sigma} + \mu_\sigma + p^n_c \right )  \charf_{t > \frac{2}{n \sigma}} 
+ e^n \charf_{0 < t < \frac{2}{n \sigma}} 
\end{aligned}
\end{equation}
with $e^n$ uniformly bounded. Accordingly, the limiting energy satisfies
\begin{equation}
\label{energybal}
E[ y(\cdot , t)  ; B] - E [ \lambda x ; B] = 
\int_0^t  \mu_{-\sigma} + \mu_\sigma + 2 (\tau_\infty - \tau(\alpha)) Y(0)  ds =:  \int_0^t  T > 0
\end{equation}
where $\mu_{-\sigma} < 0$ and $\mu_\sigma < 0$ stand for the (rate of) energy dissipation of the shocks,
$p_c = 2 ( \tau_\infty - \tau(\alpha) ) Y(0) > 0$ is the cost of opening the crack. Moreover, the total 
change of energy $T > 0$.
\end{proposition}

The only thing left to prove is the last statement. Indeed, using \eqref{mu1}, \eqref{mu2}, \eqref{kinem} and \eqref{hypo1},
we deduce
$$
\begin{aligned}
T &= \mu_{-\sigma} + \mu_\sigma + 2 (\tau_\infty - \tau(\alpha)) Y(0)
\\
&= \sigma Y(0)^2 - 2 \sigma (W(\lambda) - W(\alpha)) + 2 \tau_\infty Y(0)
\\
&= \sigma Y(0)^2 + 2 Y(0) \Big ( \tau_\infty - \frac{W(\lambda) - W(\alpha)}{\lambda - \alpha} \Big ) 
> 0 \, .
\end{aligned}
$$



\section{Cavitation}
\label{sec:cavitation}

Next, consider the case of  three or more space dimensions. The situation at hand, of a displacement field 
 discontinuous at one point, can be interpreted as cavitation.

\subsection{The setting in several space dimensions}
\label{sec:intro3d}
We begin by describing the multi--dimensional setting of the problem.
For $d\ge3,$ consider the system of elastodynamics
\begin{equation}
 \label{eq:elas3d}
{\bf y}_{tt} - \operatorname{div}({\bf S} (\nabla {\bf y}))=0 , \quad {\bf x} \in \R^d , \;  t \in \R,
\end{equation}
where ${\bf y}:\R^d \times \R \rightarrow \R^d$, subject to
the initial and boundary conditions
\begin{align}
\label{eq:ic3d}
{\bf y}({\bf x},t) = \lambda {\bf x} & \text{ for } t \leq 0 \, , \\
\label{eq:bc3d}
{\bf y}({\bf x},t) = \lambda {\bf x} & \text{ for } |{\bf x}| > \bar r t \, , 
\end{align}
for a given number $\lambda>0$ and some $\bar r >0.$ The homogeneously deformed state 
$\bar {\bf y} = \lambda { \bf x}$ is a particular solution of this problem. 

In analogy to the one-dimensional case ({\it cf.}  \eqref{elasys} ) we introduce 
the velocity ${\bf v} =  {\bf y}_t$ and the deformation gradient ${\bf F}= \nabla {\bf y}$ and observe that
\eqref{eq:elas3d} can be written as a system of conservation laws
\begin{align}\label{eq:elasys3d}
 {\bf v}_t - \operatorname{div}({\bf S}({\bf F}))&=0\\
 {\bf F}_t - \nabla {\bf v} &=0.
\end{align}
The Piola--Kirchhoff stress ${\bf S}$ is given by the assumption of hyperelasticity \eqref{eq:energy3d1},
with an energy  density $W$ that is isotropic and homogeneous. For isotropic and frame-indifferent
stored energies, 
$W({\bf F})=  \Phi (\lambda_1,\dots,\lambda_d)$ for some symmetric function $\Phi ( \lambda_1,\dots,\lambda_d )$
of the eigenvalues $\lambda_i$, $i = 1, ..., d$,  of the matrix  $\sqrt{{\bf F \bf F}^ T}$, (see \cite{TN65}).

In what follows we consider free energies $W$ having the special form \eqref{eq:energy3d2},
\begin{equation}
W({\bf F}) =  \Phi (\lambda_1,\dots,\lambda_d)
 = \frac{1}{2} \sum_{i=1}^d \lambda_i^2 + h\left(\prod_{i=1}^d \lambda_i\right)
\tag {H$_1$}
\end{equation}
where $h: \R_+ \longrightarrow \R_+$ satisfies 
\begin{equation}
 h''>0 \, ,  \quad  h'''<0 \,, \quad  \lim_{v \rightarrow 0} h(v) = \lim_{v \rightarrow \infty} h(v) = \infty  \, .
 \tag {H$_2$}
\end{equation}
The hypothesis $h'' >0$ relates to polyconvexity, while $h''' < 0$ indicates elasticity of softening type.
These hypotheses are analogous to the hypothesis  \eqref{hypo1} in
section \ref{sec-fracture} placed on the stored energy $W$ in the one-dimensional case.

We consider radially symmetric motions,
\begin{equation}
 \label{eq:radsym}
{\bf y}({\bf x},t) = w(R,t) \frac{\bf x}{R} \quad \mbox{with $R = |{\bf x}|$} \, .
\end{equation}
Then $w: \R_+ \times \R \longrightarrow \R_+$ satisfies the differential equation \eqref{eq:weq}
with the initial and boundary conditions
\begin{equation}\label{eq:wcond}
w(R,t) = \lambda R \text{ for } t \leq 0  \, , \quad  w(R,t) = \lambda R \text{ for } R >  {\bar r} t \, .
\end{equation}
Obviously, ${\bf v}= w_t(|{\bf x}|,t) \frac{\bf x}{|{\bf x}|} $ and it is derived in \cite{ball82} that
\begin{equation}\label{eq:eigenvalues}
 \lambda_1({\bf x},t) = w_R(|{\bf x}|,t),\ \lambda_2({\bf x},t)=\dots=\lambda_d({\bf x},t)= \frac{w(|{\bf x}|,t)}{|{\bf x}|}
\end{equation}
and
\begin{align}
\label{eq:symtau}
 {\bf S} ({\bf F}) &= \frac{\del \Phi}{\del \lambda_1} (\lambda_1, ... , \lambda_d )  
      \frac{\bf x}{|{\bf x}|}  \otimes \frac{\bf x}{|{\bf x}|}  
 + \frac{\del \Phi}{\del \lambda_2} (\lambda_1, ... , \lambda_d )    
  \left( I - \frac{\bf x}{|{\bf x}|}  \otimes \frac{\bf x}{|{\bf x}|} \right)
 \nonumber
 \\
 &= 
 \left(\lambda_1 + \lambda_2^{d-1} h'(\lambda_1 \lambda_2^{d-1})   \right)
\frac{\bf x}{|{\bf x}|}  \otimes \frac{\bf x}{|{\bf x}|} 
+
\left(\lambda_2 + \lambda_1 \lambda_2^{d-2} h'(\lambda_1 \lambda_2^{d-1}) \right) \left( I - \frac{\bf x}{|{\bf x}|}  \otimes \frac{\bf x}{|{\bf x}|} \right)
\end{align}
where $\lambda_1, ... , \lambda_d$ are given by \eqref{eq:eigenvalues} and $I \in \R^{d \times d}$ is the identity matrix.

In \cite{ps88, ps98}, K.A. Pericak-Spector and S. Spector  construct a self-similar weak solution of
\eqref{eq:weq}, \eqref{eq:wcond} satisfying $w(0,t)>0$ for $t>0$ and  corresponding to a cavity 
forming at time $t= 0$.
Their solution is based on the {\it ansatz} 
\begin{equation}
 \label{eq:rdef}
w(R,t)= t \, r \Big ( \frac{R}{t} \Big ) \, ,
\end{equation}
and is constructed under the framework \eqref{eq:energy3d2}, \eqref{eq:defh} in \cite{ps88}; extensions to
more general free energies are carried out in \cite{ps98}.
It is also shown in \cite{ps88,ps98} that the self-similar solution with cavity has smaller total
mechanical energy than  the associated homogeneously deformed state from where it emerges. 
The loss of energy is due  to the energy dissipation of a shock generated at the onset of cavitation 
and propagating outwards ahead of the cavity, see \cite[Thm 7.2]{ps88} .

The cavitating solution has the special property  that the Cauchy stress at the surface of the cavity vanishes 
and thereby all integrals involved in defining weak solutions and the energy of the solution are well defined.
However, the energy of this solution does not contain any contribution reflecting work which is needed to create the cavity.
Here, we will propose a different concept of solution and energy which accounts for the energetic cost for creating
the cavity. Following the ideas in section \ref{sec-fracture}, we will study 
the emergence of  the cavitating solution of \cite{ps88} as the limit of continuous approximate solutions.

We start our considerations by defining sequences of approximating functions for the class of spherically symmetric
functions in \eqref{eq:radsym}:  The function $w(R,t)$ is first extended for $R \in \R$ by using an odd extension, 
$w(-R,t) = - w(R,t)$. Consider next a symmetric mollifier $\phi \in C_c^\infty(\R) $ satisfying $\phi \geq 0,$ $\int \phi=1,$ $\operatorname{supp}(\phi)\subset [-1,1]$, and $\phi(x)=\phi(-x)$. We impose the additional restriction $\phi (0) > 0$
so that the approximate solution will detect the cavity.

Let $\phi_n( R )=n\phi(n R )$ and 
for $ w\in L^1_{loc}(\R_+ \times \R)$ define
\begin{equation}
 \label{eq:wndef}
w^n(R,t)= \int_0^\infty \phi_n(R - \tilde R) w(\tilde R,t) d\tilde R - \int_0^\infty \phi_n(R + \tilde R) w(\tilde R,t) d\tilde R.
\end{equation}
As $w$ is an odd extension,  \eqref{eq:wndef} is in fact the standard mollification of $w$  
with respect to the radial  variable $R$ .
The symmetry of $\phi$ and \eqref{eq:wndef} imply $w^n(0,t)=0.$

For a function ${\bf y} \in L_{loc}^1(\R^d \times \R,\R^d)$ of the form \eqref{eq:radsym} we define 
\begin{equation}\label{eq:yndef}
 {\bf y}^n ({\bf x},t)= w^n(|{\bf x}|,t)\frac{{\bf x}}{|{\bf x}|}.
\end{equation}
It is important to keep in mind that with this definition analogous versions of   \eqref{eq:eigenvalues} and \eqref{eq:symtau} hold.

\begin{definition}\label{def:slic3d}
Let ${\bf y} \in L^\infty_{loc}(\R; L_{loc}^1(\R^d;\R^d))$ of the form  \eqref{eq:radsym} with $w(\cdot, t)$ monotone increasing
satisfy  ${\bf y}({\bf x},t)= \lambda {\bf x}$ for $t \leq 0$ and for  $|{\bf x}| > \bar r t \, , \, t > 0$
for some $\bar r>0$. The function ${\bf y}$  is called a  singular limiting induced from continuum ({\it slic})-solution 
of \eqref{eq:elas3d} if ${\bf y}^n$ defined by \eqref{eq:yndef}, \eqref{eq:wndef} satisfies
$$
\det \nabla {\bf y}^n \geq \varepsilon_n >0 \text{ for every } n \in \mathbb{N}
$$ 
and
\[ \int_\R \int_{\R^d} {\bf y}^n \partial_{tt} {\boldsymbol \psi} + {\bf S} (\nabla {\bf y}^n) : \nabla {\boldsymbol \psi} \, d{\bf x} dt \rightarrow 0  \, , \qquad \mbox{as $n \rightarrow \infty$}, \]
holds for all $\phi \in C_c^\infty(\R)$ positive, symmetric, mollifiers with $\phi(0) > 0$, and for 
$\psi \in C_c^2(\R^d \times \R,\R^d)$.
\end{definition}

This definition takes into account the layer structure at the cavity.
In agreement with the definition of slic-solution (and of the energy in one space dimension) we define the energy of a multidimensional slic-solution.

\begin{definition}\label{def:energy3d}
The energy of a slic-solution ${\bf y} \in  W^{1,\infty}_{loc}(\R; L_{loc}^1(\R^d;\R^d))$ of the form \eqref{eq:radsym} 
in some domain $B \subset \R^d$ and for a.e. $t \in \R$ is defined as
\[ E[{\bf y},B] (t) :=  \lim_{n \to \infty} \int_B \frac{1}{2} |{\bf y}^n_t({\bf x},t)|^2 + W(\nabla {\bf y}^n({\bf x},t)) d{\bf x} \]
where ${\bf y}^n$ is defined by \eqref{eq:yndef}, \eqref{eq:wndef}.
\end{definition}

\medskip
Due to the regularity of ${\bf y}$ the quantities $\del_t {\bf y}^n$ and $\nabla {\bf y}^n$ are essentially bounded for
$(x,t)$ taking values in any compact set.
In more general situations, it might be necessary to mollify in both the space and  time variables 
the function ${\bf y}$. However, in hyperbolic problems space regularization also entails time regularization
and thus mollification in time might be at times avoided. This is certainly the case when one studies self-similar solutions.

\subsection{Solutions with cavities and their averagings }
\label{sec:prelim3d}
In the sequel, we study how the cavitating solution constructed  in \cite{ps88}
behaves with respect to Definitions \ref{def:slic3d} and \ref{def:energy3d}.
In preparation, we state some properties for the  solution from \cite{ps88}. Recall that 
$$
{\bf y} ( {\bf x}, t) = w(R,t) \, \frac{ {\bf x} }{ R } = t r \big ( \frac{R}{t} \big ) \frac{ {\bf x} }{ R } = \frac{r(s)}{s} {\bf x}\, ,
$$
with $R = |{\bf x}|$, $s = \frac{R}{t}$, and that the specific volume is given by
\begin{equation}
v := \det \nabla {\bf y} = \prod_{i=1}^d \lambda_i  = w_R \big ( \frac{w}{R} \big )^{d-1} = r'  \big ( \frac{r}{s} \big )^{d-1} .
\end{equation}

\begin{lemma}[Pericak-Spector,Spector]\label{lem:psest}
 The solution given in \cite{ps88} satifies
\begin{enumerate}
 \item $H \leq v \leq \lambda^d$ where $H$ is uniquely determined by $h'(H)=0$.
\item There is a minimal $\bar r>0$ such that ${\bf y}({\bf x},t) = \lambda {\bf x}$ for $|{\bf x}|>\bar r t$. We will denote it by $\sigma.$
\item The derivative $r'$ is non-negative and strictly increasing on $(0,\sigma]$, in particular $r'(\sigma-)< \lambda.$
\item The specific volume $v(R,t)$ is strictly increasing in $R$ for $t>0$ and $ 0 < R \leq \lambda$.
\item For any $0 \not= {\bf x} \in \R^d$ and $t>0$ it is $r(\frac{|{\bf x}|}{t})> \lambda \frac{|{\bf x}|}{t} >r'(\frac{|{\bf x}|}{t}) \frac{|{\bf x}|}{t}.$
\end{enumerate}
\end{lemma}

\noindent
For the proof of Lemma \ref{lem:psest} the reader is referred  to \cite{ps88}  and specifically to Proposition 6.2 for Assertion
2;  to Propositions 6.2 and 6.3 and (6.1) for Assertions 3 and 5; and to Proposition 6.5 (see formulas (6.19), (6.20) of \cite{ps88}) for Assertions 1 and 4.

\smallskip

Let ${\bf y}^n$ be the approximate solution defined by \eqref{eq:yndef} and \eqref{eq:wndef}, and let us fix
the following notation: For every quantity related to ${\bf y}$ the same symbol with superscript $n$ denotes 
the corresponding quantity related to ${\bf y}^n.$
 Due to \eqref{eq:yndef}, we have
\begin{equation}
\label{eq:approx}
\lambda_1^n = w^n_R \, , \quad \lambda_2^n = ... = \lambda_d^n = \frac{w^n}{R} \, , 
\quad v^n = w^n_R \Big ( \frac{w^n}{R} \Big )^{d-1} \, .
\end{equation}
We proceed to  derive uniform in $n$ bounds for the velocity, the eigenvalues and the specific volume
of the approximate solutions.

\begin{lemma}\label{lem:vest3d}
 There exist constants $c_1,c_2 >0$ such that for $R> \frac{1}{n}, \, t>0$
 \[   c_1 \leq v^n(R,t) \leq c_2.\]
\end{lemma}
\proof
The domain $R> \frac{1}{n}$, $t > 0$ can be split in two ranges:
$\sigma t + \frac{1}{n}< R$ and $ \frac{1}{n} < R < \sigma t + \frac{1}{n}$.
In the former case $v^n(R,t) =v(R,t)=\lambda^d$.

So we only consider the range $ \frac{1}{n} < R < \sigma t + \frac{1}{n}$ and 
start with proving the lower bound for $v^n$. We infer from $v(R,t) \geq H$ that
\begin{equation}\label{210601}
 w_R(R,t) \geq \frac{H R^{d-1}}{w(R,t)^{d-1}}
\end{equation}
which implies, upon using \eqref{eq:wndef}, the monotonicity property of $w$ and Jensen's inequality, that
\begin{align}\label{210602}
 w^n_R(R,t) &= \int_0^\infty \phi_n(R-\tilde R) w_R(\tilde R,t) d\tilde R
 \nonumber 
 \\
&\geq   \int_0^\infty \phi_n(R-\tilde R)  \frac{H \tilde R^{d-1}}{w(\tilde R,t)^{d-1}} 
 d\tilde R 
\nonumber 
\\
&\geq    \frac{H }{w( R + \frac{1}{n} ,t)^{d-1}}  \int_0^\infty \phi_n(R-\tilde R)  \tilde R^{d-1}  d\tilde R
\nonumber 
\\
&\geq   \frac{H }{w( R + \frac{1}{n} ,t)^{d-1}}  \left ( \int_0^\infty \phi_n(R-\tilde R)  \tilde R  d\tilde R \right )^{d-1}
\nonumber  
\\
&= \frac{H R^{d-1}}{w(R+ \frac{1}{n},t)^{d-1}} \, .
\end{align}
Moreover,
\begin{multline}\label{210603}
 w^n(R,t) = \int_0^\infty \phi_n(R - \tilde R) w(\tilde R,t) d \tilde R \\
 \geq \int_R^\infty \phi_n(R - \tilde R) w(\tilde R,t) d \tilde R
\geq \int_R^\infty \phi_n(R - \tilde R) w( R,t) d \tilde R = \frac{1}{2} w(R,t).
\end{multline}
Combining \eqref{210602} and \eqref{210603} yields
\begin{equation}\label{210604}
 v^n(R,t) \geq \frac{H R^{d-1}}{w(R+\frac{1}{n},t)^{d-1}}  \left(\frac{1}{2} \frac{w(R,t)}{R}\right)^{d-1} 
= \frac{H}{2^{d-1}} \left( \frac{w(R,t)}{w(R+\frac{1}{n},t)}\right)^{d-1}.
\end{equation}
It thus suffices to obtain a lower bound for $\frac{w(R,t)}{w(R+\frac{1}{n},t)}$.
We distinguish the cases $t \geq \frac{1}{n}$ and $t \leq \frac{1}{n}$ still relying on $\frac{1}{n} \leq R \leq \sigma t + \frac{1}{n}$,
 the ansatz \eqref{eq:rdef} and Lemma \ref{lem:psest}.
In case $t\geq \frac{1}{n}$ we find
\begin{equation}\label{210605}
 \frac{w(R,t)}{w(R+\frac{1}{n},t)} = \frac{r(\frac{R}{t})}{r(\frac{R+\frac{1}{n}}{t})} \geq \frac{r(0)}{\lambda \frac{\sigma t + \frac{2}{n}}{t}} \geq
\frac{r(0)}{\lambda (\sigma  + \frac{2}{nt})} \geq \frac{r(0)}{\lambda (\sigma  + 2 )} >0.
\end{equation}
For $t\leq\frac{1}{n}$ we obtain
\begin{equation}
 \label{210606}
\frac{w(R,t)}{w(R+\frac{1}{n},t)} = \frac{r(\frac{R}{t})}{r(\frac{R+\frac{1}{n}}{t})} \geq\frac{ \lambda \frac{R}{t}}{\lambda \frac{\sigma t + \frac{2}{n}}{t}} \geq
\frac{R}{\sigma t + \frac{2}{n}} \geq \frac{\frac{1}{n}}{\sigma \frac{1}{n}+\frac{2}{n}} \geq \frac{1}{\sigma +2} >0.
\end{equation}
Combining \eqref{210604},\eqref{210605} and \eqref{210606} we have a lower bound for $v^n.$

Next, we show an upper bound for $v^n$ always in the range $\frac{1}{n} \leq R \leq \sigma t + \frac{1}{n}$.
Observe that
$v(R,t) \leq \lambda^d$ implies
\begin{equation}
 \label{210607}  w_R(R,t) \leq \lambda^d \left( \frac{R}{w(R,t)}\right)^{d-1},
\end{equation}
and that proceeding along the lines of the derivation of \eqref{210602}, \eqref{210603} we obtain
\begin{align}
 \label{210608}  
 w^n_R(R,t) 
& \le \frac{\lambda^d }{w(R-\frac{1}{n},t)^{d-1}} \int_0^\infty \phi_n (R- \tilde R) {\tilde R}^{d-1} \, d\tilde R
 \leq \lambda^d \left( \frac{R+\frac{1}{n}}{w(R-\frac{1}{n},t)}\right)^{d-1}
 \\
  \label{210608b}  
 w^n (R,t) &= \int_0^\infty \phi_n (R- \tilde R) w(\tilde R, t) \, d\tilde R \le w(R + \frac{1}{n} , t).
\end{align}
This implies
\begin{equation} 
\label{210609}
 v^n(R,t) = w^n_R \Big ( \frac{w^n}{R} \Big )^{d-1}  \leq \lambda^d \left( \frac{R+\frac{1}{n}}{w(R-\frac{1}{n},t)} \frac{w(R+\frac{1}{n},t)}{R}\right)^{d-1}
\\ \leq \lambda^d 2^{d-1} \left( \frac{r(\frac{R+\frac{1}{n}}{t})}{r(\frac{R-\frac{1}{n}}{t})}\right)^{d-1}.
\end{equation}
For $t\geq\frac{1}{n}$ we have
\begin{equation}\label{210610}
 \frac{r(\frac{R+\frac{1}{n}}{t})}{r(\frac{R-\frac{1}{n}}{t})} \leq \frac{\lambda (\frac{\sigma t + \frac{2}{n}}{t})}{r(0)} \leq \frac{\lambda ( \sigma + \frac{2}{nt})}{r(0)} 
\leq \frac{\lambda ( \sigma + 2)}{r(0)} .
\end{equation}
For $t\leq\frac{1}{n}$ and for $\frac{1}{n} \le R \le \sigma t + \frac{1}{n}$ we estimate $v^n$ without using \eqref{210609}.
  We rely on $w^n_R(R,t) \leq \lambda$ and
\begin{equation}\label{210611}
 \frac{w^n(R,t)}{R} \leq \frac{w (R+ \frac{1}{n},t)}{R} \leq \frac{w(0,t) + \lambda(R + \frac{1}{n})}{R} \leq \frac{t}{R} r(0) + \lambda(1+ \frac{1}{nR}) \leq r(0)+2\lambda.
\end{equation}
Combining \eqref{210609}, \eqref{210610} and \eqref{210611} we see that $v^n$ is bounded from above.
\qed


The next Lemma establishes a global bound for ${\bf y}^n_t$ and bounds for the eigenvalues outside a neighborhood of the cavity.

\begin{lemma}\label{lem:bulkest3d}
 There exists a constant $C>0$ (independent of $n$) such that
\begin{enumerate}
 \item $|{\bf y}^n_t({\bf x},t)| = |w^n_t(|{\bf x}|,t)| \leq C $ for all $0 \not= {\bf x} \in \R^d$ and $t>0.$
 \item $0 \leq \lambda_1^n({\bf x},t) \leq C $ for $|{\bf x}| > \frac{1}{n}$ and $t>0.$
 \item $0 \leq \lambda_2^n({\bf x},t) \leq C\left(1 + \frac{t}{|{\bf x|}} \right) $ for $|{\bf x}| > \frac{1}{n}$ and $t>0.$
\end{enumerate}

\end{lemma}
\proof
Using \eqref{eq:wndef} we write
\begin{align*}
w^n_t(R,t) &=\int_0^\infty \phi_n(R-\tilde R) w_t(\tilde R,t) d \tilde R- \int_0^\infty \phi_n(R+\tilde R) w_t(\tilde R,t) d \tilde R .
\end{align*}
Also, from Lemma \ref{lem:psest} and \eqref{eq:rdef} we have $w_t (R, t) = r(s) - s r^\prime (s) > 0$,
$\frac{d}{ds} \big ( r - s r^\prime) = - s r^{\prime \prime} < 0$, and hence $0 < w_t (R,t) < r(0)$ for $R > 0 \, , \; t > 0$.
To show the first assertion we distinguish two cases.
For $R> \sigma t + \frac{1}{n}$ we have
\[ w^n_t(R,t) =w_t(R,t)=\partial_t(\lambda R)=0.\]
For $R< \sigma t + \frac{1}{n}$ we find
\begin{align*}
w^n_t(R,t) \le \int_0^\infty \phi_n(R-\tilde R) w_t(\tilde R,t) d \tilde R
\leq 
 \int_0^\infty \phi_n(R-\tilde R) r(0)d \tilde R 
 \leq r(0).
\end{align*}

To prove the second assertion, we note that for $R>\frac{1}{n}$
it holds
\[ w^n_R(R,t)= \int_\R \phi_n(R-\tilde R) w_R(\tilde R,t) d \tilde R\]
such that upper and lower bounds for $w_R$ are upper and lower bounds for $\lambda_1^n$.

It remains to prove the third assertion. We have
\begin{multline}
 \frac{w^n(R,t)}{R} = \frac{1}{R} \int_\R \phi_n(R - \tilde R) w(\tilde R,t) d \tilde R \leq \frac{w(R+ \frac{1}{n} , t)}{R}\\
\leq \frac{w(0,t) + \lambda (R + \frac{1}{n})}{R} = \frac{r(0)t}{R} + \lambda + \frac{\lambda}{Rn} \leq \frac{r(0)t}{R} + 2\lambda.
\end{multline}

\qed

The next lemma investigates the behaviour of eigenvalues in the vicinity of the cavity. In preparation, we define
\begin{equation}
 \label{Phidef}
\begin{aligned}
\Phi(R) &:= \frac{1}{2} \int_0^\infty \phi(R-\tilde R) - \phi(R+\tilde R) d \tilde R
\\
&= \int_0^R \phi (x) dx
\end{aligned}
\end{equation}


\begin{lemma}\label{lem:cavest3d}
 For $R < \frac{1}{n}$ and $t>0$ the following holds
\begin{enumerate}
\item $0 \le 2\phi_n(R) w(0,t) \leq  \del_R w^n(R,t) \leq 2\phi_n(R) w(0,t) + \lambda.$ 
\item $2 \frac{\Phi(nR)}{R} w(0,t)  \leq \frac{w^n(R,t)}{R} \leq 2 \frac{\Phi(nR)}{R} w(0,t) + \lambda $.
\item If $\phi(0) \ne 0$ then there exist $\eps \, , \; \delta > 0$ such that  for $| {\bf x} | < \frac{\eps}{n}$
 $$
 v^n ( {\bf x}, t) \ge 2^d n^d \delta^d w(0,t)^d \, .
 $$ 
\item  If $\phi(0) \ne 0$ there exist constants $c_1,c_2 >0$ such that 
\begin{equation}
\label{vnbound}
c_1 \leq v^n({\bf x},t) \leq c_2 (1+ t^dn^d) \, , \quad \mbox{for $|{\bf x}|< \frac{1}{n}$, $t > 0$}.
\end{equation}
\end{enumerate}

\end{lemma}
\proof
Upon differentiating \eqref{eq:wndef} and using the integration by parts formula, we obtain
\begin{equation}
\label{101}
\begin{aligned}
\del_R w^n(R,t)= 2 \phi_n(R) w(0,t) + \int_0^\infty \big [ \phi_n(R - \tilde R) + \phi_n(R + \tilde R) \big ]  w_R(\tilde R,t) d \tilde R.
\end{aligned}
\end{equation}
Recall that $0\le  \del_R w \le \lambda$. Then
\begin{align*}
0 &\le \int_0^\infty \big [ \phi_n(R - \tilde R) + \phi_n(R + \tilde R) \big ]  w_R(\tilde R,t) d \tilde R
\\
&\le \lambda \int_0^\infty  \phi_n(R - \tilde R) + \phi_n(R + \tilde R) d \tilde R
\\
&= \lambda
\end{align*}
and \eqref{101} shows Assertion 1. By \eqref{Phidef}, it also implies
\begin{align*}
\frac{\del}{\del R} \big ( 2 \Phi(nR) w(0,t) \big ) \le \del_R w^n \le \frac{\del}{\del R} \big ( 2 \Phi(nR) w(0,t) \big ) + \lambda
\end{align*}
which once integrating and noting that $w^n (0,t) = 0$ proves Assertion 2.

We turn now to Assertions 3 and 4. Observe that, since $w^n (0,t) = 0$,
\begin{equation}
\label{103}
v^n (R,t) = w^n_R (R,t) \left ( \frac{w^n(R,t)}{R} \right )^{d-1} =
w_R^n (R,t) \left ( \frac{1}{R}   \int_0^R w^n_R (\bar R, t) \, d\bar R \right )^{d-1}
\ge \inf_{0 < \bar R < R}  ( w_R^n (\bar R , t) )^d.
\end{equation}
If $\phi(0) \ne 0$ there are $\eps > 0$ and $\delta > 0$ such that $\phi (x) > \delta$ for $|x | < \eps$.
We restrict our computations to the range $R = | {\bf x} | < \frac{\eps}{n}$ and use \eqref{103} and Assertion 1 to infer
\begin{align*}
v^n (R,t) \ge \inf_{0 < \bar R < \frac{\eps}{n}}  ( w_R^n (\bar R , t) )^d
\ge 2^d \inf_{0 < \bar R < \frac{\eps}{n}} ( n \phi( n {\bar R} ) w(0,t) )^d 
\ge 2^d n^d \delta^d w(0,t)^d
\end{align*}
which shows Assertion 3.

Turning now to \eqref{vnbound}, we note that the upper bound follows from Assertions 1. and 2. via
$$
\begin{aligned}
v^n &= w^n_R \left ( \frac{w^n}{R} \right )^{d-1} 
\\
&\le \big ( 2 n \phi(nR) w(0,t) + \lambda \big ) \Big (  w(0,t) \frac{2}{R} \int_0^{nR} \phi(x) dx  + \lambda \Big )^{d-1}
\\
&\le C (1 + n^d t^d ).
\end{aligned}
$$
The lower bound is more delicate and we distinguish three cases:
\begin{itemize}
\item[(i)]  $0 < R < \frac{\eps}{n}$, $t \ge \frac{\eta}{2 n \sigma}$ for some $0< \eta \le 1$ to be selected below;
\item[(ii)]  $0 < R < \frac{\eps}{n}$, $0 <  t \le \frac{\eta}{2 n \sigma}$;
\item[(iii)] $\frac{\eps}{n} < R < \frac{1}{n}$, $t > 0$.
\end{itemize}
For $(R, t)$ in the range (i), we have
$$
w^n_R (R,t) \ge 2 \phi_n (R) w(0,t) = 2 n \phi (n R)  w(0,t) \ge 2 n t \delta r(0) \ge \frac{\eta}{\sigma} \delta r(0) > 0
$$
and the desired bound follows from \eqref{103}.

Next, let $(R, t)$ be in the range (ii) and define $\hat R = n R$ so that $0 < \hat R < \eps$. Using \eqref{101}
we obtain
\begin{align*}
w^n_R (R, t)  &\ge \int_0^\infty \big ( n \phi ( \hat R - n \bar R ) + n \phi ( \hat R + n \bar R) \big ) w_R ( \bar R , t) \, d\bar R
\\
&=  \int_0^\infty \big (  \phi ( \hat R -  \bar R ) +  \phi ( \hat R +  \bar R) \big ) w_R ( \frac{\bar R}{n} , t) \, d\bar R
\\
&\ge \int_{n \sigma t}^\infty (  \phi ( \hat R -  \bar R ) +  \phi ( \hat R +  \bar R) \big ) w_R ( \frac{\bar R}{n} , t) \, d\bar R
\\
&=  \int_{n \sigma t}^\infty (  \phi ( \hat R -  \bar R ) +  \phi ( \hat R +  \bar R) \big ) \lambda \, d\bar R
\\
&\ge  \lambda \int_{\eta/2}^\infty \phi ( \hat R -  \bar R )  \, d\bar R
\\
&= \lambda \int_{-\infty}^{\hat R - \frac{\eta}{2}} \phi (x) dx
>  \lambda \int_{-\infty}^{ - \frac{\eta}{2}} \phi (x) dx.
\end{align*}
By selecting $\eta$ sufficiently small, we ensure that $\int_{-\infty}^{ - \frac{\eta}{2}} \phi (x) dx \ge c > 0$ 
and the lower bound follows in range (ii) follows from \eqref{103}.

We next consider $(R,t)$ in range (iii). Using again \eqref{101} we have
\begin{align*}
w^n_R (R, t)  &\ge \int_0^\infty \phi_n (R - \bar R) w_R (\bar R, t) \, d\bar R 
\ge \int_{\frac{\eps}{2n}}^{1+ \frac{\eps}{2n}} \phi_n (R - \bar R) w_R (\bar R, t) \, d\bar R 
\\
&\ge  w_R ( \frac{\eps}{2n}, t) \int_{\frac{\eps}{2n}}^{1+ \frac{\eps}{2n}} \phi_n (R - \bar R) \, d\bar R 
\ge  w_R ( \frac{\eps}{2n}, t)   \int_{ - \frac{\eps}{2n}}^{  \frac{\eps}{2n}}  \phi_n (x) dx 
\\
&\ge \eps \delta \, w_R ( \frac{\eps}{2n}, t). 
\end{align*}
Moreover, using Assertion 1., \eqref{eq:wndef}, the symmetry of $\phi$ and the monotonicity of $w(\cdot, t)$,
we obtain
\begin{align*}
w^n (R,t) &\ge w^n ( \frac{\eps}{n}, t) 
    = \int_0^\infty \big ( \phi_n ( \frac{\eps}{n} - \bar R) - \phi_n ( \frac{\eps}{n} + \bar R ) \big ) w ( \bar R , t) \, d\bar R
\\
&= \int_{-\frac{ \eps}{n}}^\infty \phi_n ( -\tilde R) w \big ( \tilde R + \frac{\eps}{n} , t \big ) \, d\tilde R
    -  \int_{\frac{ \eps}{n}}^\infty \phi_n ( \tilde R) w\big ( \tilde R - \frac{\eps}{n} , t \big )  \, d\tilde R 
\\
&= \int_{-\frac{ \eps}{n}}^{\frac{1}{n}} \phi_n ( \tilde R) w \big ( \tilde R + \frac{\eps}{n} , t \big ) \, d\tilde R
    -  \int_{\frac{ \eps}{n}}^{\frac{1}{n}} \phi_n ( \tilde R) w\big ( \tilde R - \frac{\eps}{n} , t \big )  \, d\tilde R 
\\
&= \int_{-\frac{ \eps}{n}}^{\frac{ \eps}{n}} \phi_n ( \tilde R) w \big ( \tilde R + \frac{\eps}{n} , t \big ) \, d\tilde R 
+ \int_{\frac{ \eps}{n}}^{\frac{1}{n}} \phi_n ( \tilde R)
 \Big (  w \big ( \tilde R + \frac{\eps}{n} , t \big ) - w\big ( \tilde R - \frac{\eps}{n} , t \big )  \Big )  \, d\tilde R  
 \\
 & \ge \int_0^{\frac{ \eps}{n}} \phi_n ( \tilde R) w \big ( \tilde R + \frac{\eps}{n} , t \big ) \, d\tilde R 
  \\
 &\ge w \big ( \frac{\eps}{n} , t \big ) \int_0^{\frac{ \eps}{n}} \phi_n ( \tilde R)  \, d\tilde R 
 \ge 
 \eps \delta w \big ( \frac{\eps}{n} , t \big )
 \ge 
  \eps \delta w \big ( \frac{\eps}{2n} , t \big ).
\end{align*}
In turn, always in the range $\frac{\eps}{n} < R < \frac{1}{n}$,
\begin{multline}
v^n = w^n_R \left ( \frac{w^n}{R} \right )^{d-1}  
\ge (\eps \delta)^d  \frac{ w_R \big ( \frac{\eps}{2n} , t \big ) \Big(w \big ( \frac{\eps}{2n} , t \big ) \Big)^{d-1} }{R^{d-1}}
\\
\ge (\eps \delta)^d \Big ( \frac{\eps}{2} \Big )^{d-1}
 \frac{ w_R \big ( \frac{\eps}{2n} , t \big ) \Big(w \big ( \frac{\eps}{2n} , t \big )\Big)^{d-1}}{\big ( \frac{\eps}{2n}\big )^{d-1}}
= (\eps \delta)^d \Big ( \frac{\eps}{2} \Big )^{d-1} v \big ( \frac{\eps}{2n} , t \big ) 
\ge (\eps \delta)^d \Big ( \frac{\eps}{2} \Big )^{d-1} H
\end{multline}
which concludes the proof of Assertion 4.
\qed

\begin{remark} \rm
The hypothesis $\phi(0) \ne 0$ is necessary  to obtain the lower bounds stated in Assertions 3. and 4.
Indeed, when the mollifying kernel satisfies $\phi (0) = 0$ then we will show that the approximate solution
satisfies
\begin{equation}
\label{rem101}
\lim_{n \to \infty} \sup_{ |{\bf x}| < \frac{1}{n^2} }  v^n ( | {\bf x}| , t) = 0
\end{equation}
for every $t>0$. This entails that the approximate solution will have regions in space with behavior of large compression 
intermixed with neighboring regions with behavior of large extensions, which is quite unexpected for
a large extension situation as appears in cavitation. The assumption $\phi(0) \ne 0$ rules out such behavior
for the approximants.

To show \eqref{rem101}, note that due to the symmetry $\phi(-R)=\phi(R)$ the condition $\phi(0) = 0$ implies
that $\phi(R) \le c R^2$ for $R$ sufficiently small. Consider now the region $0 < R < \frac{1}{n^2}$ and note that
\eqref{101} implies
\begin{align}
w^n_R (R,t) &= 2 n \phi(n R) w(0,t) + \int_0^\infty \big [ \phi_n(R - \tilde R) + \phi_n(R + \tilde R) \big ]  w_R(\tilde R,t) d \tilde R
\nonumber
\\
&\le \frac{2 c}{n} w(0,t) + \sup_{0 < \tilde R \le \frac{1}{n}+ \frac{1}{n^2}} w_R (\tilde R, t)
\label{105}
\\
&= \frac{2 c}{n} w(0,t) + c_n (t) \longrightarrow 0 \qquad  \mbox{as $n \to \infty$}
\nonumber
\end{align}
where we have set $c_n (t) := \sup_{0 < \tilde R \le \frac{1}{n}+ \frac{1}{n^2}} w_R (\tilde R, t)$ and used 
that $w_R (0,t) = 0$. Integrating \eqref{105} gives
$$
w^n (R, t) \le \frac{2 c}{n} w(0,t) R  + c_n (t) R
$$
and thus
$$
v^n (R,t) = w^n_R \left ( \frac{w^n}{R} \right )^{d-1}  \le \left (  \frac{2 c}{n} w(0,t) + c_n (t) \right )^d
\, , \qquad \mbox{for $0 < R \le \frac{1}{n^2}$,}
$$
which shows \eqref{rem101}.

\end{remark}

\subsection{Existence of multi--dimensional slic-solutions}

In this section we investigate under which conditions the weak solutions constructed in \cite{ps88} and 
exhibiting cavitation are also a slic-solution according to Definion \ref{def:slic3d}.
This turns out to depend on the growth behaviour of $h(v)$ for $v \rightarrow \infty$.

\begin{theorem}\label{lem:slic3d}
 The weak solutions constructed in \cite{ps88} -- extended by ${\bf y}({\bf x},t)=\lambda {\bf x}$ for $t<0$ -- are slic-solutions provided
\[ \lim_{v \rightarrow \infty} \frac{h'(v^d)}{v} =0.\]
They are not slic-solutions in case
\[ \liminf_{v \rightarrow \infty} \frac{h'(v^d)}{v} >0.\]
\end{theorem}

\proof
 By construction the solutions from \cite{ps88} satisfy
\begin{equation}
 \int_0^\infty \int_{\R^d} \partial_t {\bf y} \partial_t {\boldsymbol \psi} - {\bf S}(\nabla{\bf y}) : \nabla {\boldsymbol \psi} d{\bf x} dt =0
\end{equation}
for all ${\boldsymbol \psi} \in C_c^\infty(\R^d \times [0,\infty),\R^d)$ and
\begin{equation}
 \lim_{t \rightarrow 0+} {\bf y}({\bf x},t) = \lambda {\bf x} \quad \text{and} \quad \lim_{t \rightarrow 0+}\partial_t  {\bf y}({\bf x},t) ={\bf 0}.
\end{equation}
The cavitating solution is extended to negative values of $t$ by setting ${\bf y}({\bf x},t) = \lambda {\bf x}$ for $t \le 0$. 
It is easy to check that the extended function (still denoted by) ${\bf y}$ is a weak solution on $\R^d \times \R$ and 
satisfies the weak form
\begin{equation}
  \int_\R \int_{\R^d} {\bf y} \partial_{tt} {\boldsymbol \psi} + {\bf S}(\nabla{\bf y}) : \nabla {\boldsymbol \psi} d{\bf x} dt =0
\end{equation}
for any ${\boldsymbol \psi} \in C_c^\infty( \R^d \times \R , \R^d)$.

To investigate whether ${\bf y}$ is a slic-solution it suffices to study whether
\begin{equation}\label{120616}
 \lim_{n \rightarrow \infty}   \int_\R \int_{\R^d} {\bf y}^n \partial_{tt} {\boldsymbol \psi} + {\bf S}(\nabla{\bf y}^n) : \nabla {\boldsymbol \psi} d{\bf x} dt =
  \int_\R \int_{\R^d} {\bf y} \partial_{tt} {\boldsymbol \psi} + {\bf S}(\nabla{\bf y}) : \nabla {\boldsymbol \psi} d{\bf x} dt.
\end{equation}
To check this, we decompose the integral on the left into three parts
\begin{multline}
 \int_\R \int_{\R^d} {\bf y}^n \partial_{tt} {\boldsymbol \psi} + {\bf S}(\nabla{\bf y}^n) : \nabla {\boldsymbol \psi} d{\bf x} dt
= \int_\R \int_{\R^d} {\bf y}^n \partial_{tt} {\boldsymbol \psi}  d{\bf x} dt\\
+ \int_\R \int_{\R^d}  {\bf S}(\nabla{\bf y}^n) : \nabla {\boldsymbol \psi} \charf_{\{|{\bf x}| > \frac{1}{n}\}} d{\bf x} dt
+ \int_\R \int_{\R^d}  {\bf S}(\nabla{\bf y}^n) : \nabla {\boldsymbol \psi} \charf_{\{|{\bf x}| < \frac{1}{n}\}} d{\bf x} dt
=: J_1^n + J_2^n + J^n_3.
\end{multline}
By standard properties of convolution ${\bf y}^n$ converges to ${\bf y}$ in $L^1(\operatorname{supp}({\boldsymbol \psi})).$ Hence, $J_1^n$ converges to
\[ \int_\R \int_{\R^d} {\bf y} \partial_{tt} {\boldsymbol \psi}  d{\bf x} dt.\]
We will now show that \[ J^n_2\rightarrow \int_\R \int_{\R^d}  {\bf S}(\nabla{\bf y}) : \nabla {\boldsymbol \psi}  d{\bf x} dt.\] 
The cavitating solution  $w$ constructed in \cite{ps88} is an element of $W^{1,1}((0,\infty)\times (0,\infty),\R)$
and thus $w^n \rightarrow w$ and $w^n_R \rightarrow w_R$ pointwise almost everywhere; accordingly 
$\nabla {\bf y}^n \rightarrow \nabla {\bf y}$ pointwise almost everywhere.
It remains to show that the integrand in $J^n_2$ is uniformly bounded by an $L^1$-function. To shorten notation 
we define $A:= \supp({\boldsymbol \psi})$ and compute
\begin{equation}
\begin{split}
 &|{\bf S}(\nabla {\bf y}^n) : \nabla {\boldsymbol \psi} \charf_{\{|{\bf x}| > \frac{1}{n}\}}| 
 \\
&\leq  \|\nabla {\boldsymbol \psi} \|_\infty \charf_A
\left(\lambda_1^n(|{\bf x}|,t)+ \lambda_2^n(|{\bf x}|,t) 
+\big[ \lambda_1^n(|{\bf x}|,t) \lambda_2^n(|{\bf x}|,t)^{d-2} + \lambda_2^n(|{\bf x}|,t)^{d-1}\big] h'(v^n(|{\bf x}|,t)) \right)
\charf_{\{|{\bf x}| > \frac{1}{n}\}}
\\
&\stackrel{\text{Lemma } \ref{lem:vest3d}}{\leq} C\|\nabla {\boldsymbol \psi} \|_\infty \charf_A
\left(\lambda_1^n(|{\bf x}|,t)+ \lambda_2^n(|{\bf x}|,t) 
+ \lambda_1^n(|{\bf x}|,t) \lambda_2^n(|{\bf x}|,t)^{d-2} + \lambda_2^n(|{\bf x}|,t)^{d-1} \right)
\charf_{\{|{\bf x}| > \frac{1}{n}\}}
\\
&\stackrel{\text{Lemma } \ref{lem:bulkest3d}}{\leq} C\|\nabla {\boldsymbol \psi} \|_\infty \charf_A
\max\left\{  1, 1+ \frac{t}{R} , \left(  1+ \frac{t}{R}\right)^{d-2} , \left( 1+ \frac{t}{R}\right)^{d-1}
\right\}
\end{split}
\end{equation}
which is in $L^1.$

Hence, \eqref{120616} holds if and only if $J^n_3 \rightarrow 0$.
Let us start with the case $\lim_{v \rightarrow \infty} \frac{h'(v^d)}{v} =0$. Then, it is sufficient to show that for every $T>0$
\begin{equation}\label{210617}
 \lim_{n \rightarrow \infty} \int_{-T}^T \int_{\{ |{\bf x}|< \frac{1}{n}\}} | {\bf S} (\nabla {\bf y}^n)| d{\bf x} dt =0. 
\end{equation}
As ${\bf S} (\nabla {\bf y})$ is diagonalisable this holds provided
\begin{align}
 \label{210618}
\lim_{n \rightarrow \infty} \int_{-T}^T \int_{\{ |{\bf x}|< \frac{1}{n}\}} \lambda_1^n(|{\bf x}|,t) d{\bf x} dt&=0 \\
 \label{210619}\lim_{n \rightarrow \infty} \int_{-T}^T \int_{\{ |{\bf x}|< \frac{1}{n}\}} \lambda_2^n(|{\bf x}|,t) d{\bf x} dt&=0 \\
 \label{210620}\lim_{n \rightarrow \infty} \int_{-T}^T \int_{\{ |{\bf x}|< \frac{1}{n}\}} \lambda_1^n(|{\bf x}|,t)\lambda_2^n(|{\bf x}|,t)^{d-2} h'(v^n(|{\bf x}|,t)) d{\bf x} dt&=0 \\
 \label{210621}\lim_{n \rightarrow \infty} \int_{-T}^T \int_{\{ |{\bf x}|< \frac{1}{n}\}} \lambda_2^n(|{\bf x}|,t)^{d-1} h'(v^n(|{\bf x}|,t)) d{\bf x} dt&=0 . 
\end{align}
In \eqref{210618} and \eqref{210619} the integrand is bounded by $Cn$ according to Lemma \ref{lem:cavest3d} and the volume of the integrational domain is proportional to $n^{-d}.$
For the fourth limit we find
\begin{equation}
 \begin{split}
& | \lim_{n \rightarrow \infty} \int_{-T}^T \int_{\{ |{\bf x}|< \frac{1}{n}\}} \lambda_2^n(|{\bf x}|,t)^{d-1} h'(v^n(|{\bf x}|,t)) d{\bf x} dt|\\ 
&\leq \lim_{n \rightarrow \infty} \int_{-T}^T \int_0^{\frac{1}{n}} (Cn)^{d-1} \max\{ C, h'(c(1+t^dn^d))\} R^{d-1}  dR dt\\
& \leq \lim_{n \rightarrow \infty} \int_{-T}^T \int_0^1 C^{d-1} \max\{ C, h'(c(1+t^dn^d))\} \frac{1}{n} R^{d-1}  dR dt\\
& \leq C  \int_{-T}^T \int_0^1 \lim_{n \rightarrow \infty}  \frac{h'(c(1+t^dn^d))}{\sqrt[d]{c(1+t^dn^d)}} \frac{\sqrt[d]{c(1+t^dn^d)}}{n}R^{d-1}  dR dt=0,
 \end{split}
\end{equation}
by the dominated convergence theorem.
The third term can be treated analogously.
Thus, the first assertion of the Lemma is proven.

In case $\liminf_{v \rightarrow \infty} \frac{h'(v^d)}{v} >0$
we consider special test functions ${\boldsymbol \psi}$ which satisfy ${\boldsymbol \psi}({\bf x},t)= \zeta(t){\bf x}$ for $|{\bf x}| <1$
with $\supp{\zeta} \subset [0,1]$ and $\zeta \geq 0 .$
For such a test function
\begin{equation}\label{210622}\begin{split}
J_3^n &=  \int_{-T}^T \int_{\{ |{\bf x}|< \frac{1}{n}\}}  {\bf S} (\nabla {\bf y}^n) : \nabla {\boldsymbol \psi} d {\bf x} dt
 =  \int_{-T}^T \int_{\{ |{\bf x}|< \frac{1}{n}\}}  {\bf S} (\nabla {\bf y}^n) : I  \zeta(t) d {\bf x} dt\\
&\geq  \int_{-T}^T \int_{\{ |{\bf x}|< \frac{1}{n}\}}  \lambda_2^n(|{\bf x}|,t)^{d-1} 
h'(v^n(|{\bf x}|,t)) \zeta(t)  d {\bf x} dt\\
&\qquad + \int_{-T}^T \int_{\{ |{\bf x}|< \frac{1}{n}\}} (d-1) \lambda_1^n(|{\bf x}|,t) \lambda_2^n(|{\bf x}|,t)^{d-2} 
h'(v^n(|{\bf x}|,t)) \zeta(t)  d {\bf x} dt 
\\
&=: I_1^n + I_2^n
\end{split}\end{equation}
making use of \eqref{eq:symtau} and noting that all eigenvalues are positive.
Let us replace $h'$ by a function which is easier to handle.
In case $A : = \liminf_{v \rightarrow \infty} \frac{h'(v^d)}{v} < \infty$ we define
\[ g(v):= \min\left\{ h'(v),  A  \sqrt[d]{v}\right\}.\]
For $\liminf_{v \rightarrow \infty} \frac{h'(v^d)}{v} = \infty$ we employ the definition
\[ g(v) := \left\{\begin{array}{ll}
                   h'(v) & \text{for } v \leq \bar v\\
                 c \sqrt[d]{v} & \text{for } v > \bar v
                  \end{array}\right.
\]
where $\bar v,c >0$ are chosen such that $h'(v)> c \sqrt[d]{v} $ for $v > \bar v.$
In both cases $g$ has the following properties, which we will exploit in the sequel
\begin{equation}\label{210623}
\begin{split}
&g(v) \leq h'(v) \quad \forall \, v>0, 
\\
&\lim_{v \rightarrow \infty} \frac{g(v^d)}{v} \text{ exists and is finite and positive},\\
& \exists C>0 \text{ such that }
|g(v^d)| \leq C v \ \forall v \in [ \min_{n,{\bf x},t} v^n({\bf x},t), \infty).
\end{split}
\end{equation}
Combining \eqref{210622} and \eqref{210623} we get
\begin{equation}\label{210624}
I_1^n
\geq  \omega_d \int_{-T}^T \int_0^1  \lambda_2^n(\frac{R}{n},t)^{d-1}
g(v^n(\frac{R}{n},t)) R^{d-1} n^{-d} d R \; \zeta(t) dt
\end{equation}
where $\omega_d$ is the surface area of the unit sphere.
Let us note that the absolute value of the integrand on the right hand side of \eqref{210624} is bounded by
\[ (Cn)^{d-1}
 R^{d-1} n^{-d} \max\{ C, g(Cn^d)\} \leq C.\]
Therefore integral and limit commute and we obtain
 \begin{equation}\label{210625}\begin{split}
\liminf_{n \rightarrow \infty} I_1^n
&\geq  \omega_d \int_{-T}^T \int_0^1 \lim_{n \rightarrow \infty}  \lambda_2^n(\frac{R}{n},t)^{d-1} 
g(v^n(\frac{R}{n},t)) R^{d-1} n^{-d} d R \zeta(t) dt\\
&= \omega_d \int_{-T}^T \int_0^1 \lim_{n \rightarrow \infty} \Big[
\left( 2n \frac{\Phi(R)}{R} w(0,t) + \mathcal{O}(1)\right)^{d-1}
 \frac{R^{d-1}}{n^{d-1}}  \\
&\qquad \qquad \frac{g\left(  (2 n \phi(R) w(0,t) + \mathcal{O}(1)) \left( 2n \frac{\Phi(R)}{R} w(0,t) + \mathcal{O}(1)\right)^{d-1} \right) }{n} \Big] dR \zeta(t)dt\\
&=  \omega_d \int_{-T}^T \int_0^1 
2^{d-1}  w(0,t)^{d-1}  \left(\frac{\Phi(R)}{R}\right)^{d-1}
  R^{d-1} \\
&\qquad\qquad\qquad \qquad
 \lim_{\hat v \rightarrow \infty} \frac{g(\hat v^d)}{\hat v}  2 w(0,t)\sqrt[d]{\phi(R) \left( \frac{\Phi(R)}{R}\right)^{d-1}} dR \zeta(t) dt
 \; > \; 0 
\end{split} \end{equation}
for a suitable choice of $\zeta.$
By an analogous calculation we derive a lower bound for the term $I_2^n$ in \eqref{210622},
\begin{multline}
\liminf_{n \rightarrow \infty} I_2^n
\ge  \omega_d \int_{-T}^T \int_0^1 
2^{d-1}  w(0,t)^{d-1}  (d-1)  \phi(R) \left(\frac{\Phi(R)}{R}\right)^{d-2}
 R^{d-1} \\
\lim_{\hat v \rightarrow \infty} \frac{g(\hat v^d)}{\hat v}  2 w(0,t)\sqrt[d]{\phi(R) \left( \frac{\Phi(R)}{R}\right)^{d-1}} dR \zeta(t) dt
\; > \; 0.
\end{multline}

\subsection{The energy needed to create a cavity}
\label{sec:energy3d}

We will now study the energy of the solutions constructed in \cite{ps88} using our notion of energy, cf. Definition \ref{def:energy3d}.
We will only consider energies in domains $B$ with finite volume.

Let us note first that if $B$ is a set with Lipschitz boundary containing a neighborhood of the whole wave fan at time t, then
\begin{equation}
 \lim_{n \rightarrow \infty} \int_B \operatorname{div}({\bf S}(\nabla {\bf y}^n) \partial_t {\bf y}^n ) d {\bf x} = 
 \lim_{n \rightarrow \infty}\int_{\partial B} {\bf n}^T {\bf S}(\nabla {\bf y}^n) \partial_t {\bf y}^n =0
\end{equation}
where ${\bf n}$ is the unit outer normal to $\partial B$.
Thus, in the limit there is no work performed by exterior forces.

Our first result shows that, if the  material is not sufficiently weak, cavitating solutions  have infinite energy for arbirarily small positive times, although their energy in $B$ at time $t=0$ is
given by $|B| W(\lambda)$, i.e. is finite.

\begin{proposition}\label{lem:energy3d1}
If $\lim_{v \rightarrow \infty}\frac{h(v)}{v} = \infty$ the energy of the cavitating solution constructed in \cite{ps88} 
in the sense of the Definition \ref{def:energy3d} satisfies
\begin{equation}
 E[{\bf y},B](t)= \infty \quad \text{ for every } t>0.
\end{equation}
\end{proposition}

\proof
For approximate solutions ${\bf y}^n$  generated by mollifiers with $\phi(0) \ne 0$, 
we fix $\eps$ and $\delta$ as in Lemma \ref{lem:cavest3d} and recall that $v^n$ satisfies Assertions 3. and 4. of that lemma.
As $B$ contains the whole wave fan at time $t$, it also contains $B(n):=B_{\frac{\eps}{n}}(0)$ for $n$ small enough. The 
non--negativity of $W$ implies
\begin{equation}\label{210626}
 E[{\bf y},B](t) \geq \lim_{n \rightarrow \infty} \int_{B(n)} h(v^n(|{\bf x}|,t)) d {\bf x} 
= \lim_{n \rightarrow \infty} \omega_d \int_0^{\eps} h(v^n(\frac{R}{n},t)) \frac{R^{d-1}}{n^d} d R.
\end{equation}
Moreover, for $0< R < \eps$, we have
\[ v^n(\frac{R}{n},t) \geq (2\delta n)^d w(0,t)^d. \]
By Fatou's lemma
\begin{align*}
  E[{\bf y},B](t)  &\geq 
  \liminf_{n \rightarrow \infty} \omega_d \int_0^{\eps} h \big((2\delta n )^d w(0,t)^d \big) \frac{R^{d-1}}{n^d} d R
  \\
&\ge \omega_d \int_0^{\eps} \liminf_{n \rightarrow \infty} \frac{h \big((2\delta n )^d w(0,t)^d \big)}{n^d} R^{d-1} d R = \infty.
\end{align*}

\qed

For weaker materials, i.e. $\lim_{v \rightarrow \infty} \frac{h(v)}{v}< \infty$, the situation changes and the energy of cavitating solutions becomes finite.
\begin{proposition}\label{lem:energy3d2}
 Let $L:= \lim_{v \rightarrow \infty} \frac{h(v)}{v} $ be finite and let $B=B_{\bar R}(0)$ for some $\bar R>0$ contain a neighborhood of the whole wave fan at time $t$. Then, the energy of 
the weak solution from \cite{ps88} - in the sense of Definition \ref{def:energy3d} - satisfies
\begin{equation}
 E[{\bf y},B] (t)= E[\lambda {\bf x},B](t) + \frac{t^d \sigma^d \omega_d}{d} J + \frac{t^d \omega_d}{d} r(0)^d L,
\end{equation}
where 
\begin{multline}
 J:= \frac{1}{2} w_R(t\sigma-,t)^2 + h(w_R(t\sigma-,t)\lambda^{d-1}) - \frac{1}{2} \lambda^2 - h(\lambda^d)\\
  + \frac{1}{2} \left[   w_R(t \sigma-,t) + h'(w_R(t \sigma-,t) \lambda^{d-1} ) \lambda^{d-1}+ \lambda + h'(\lambda^d)\lambda^{d-1}   
\right]  
(\lambda -  w_R(t \sigma-,t) )
\end{multline}
is the energy dissipation of the outgoing shock.

\end{proposition}

\proof
 Taking into account Proposition 7.1 from \cite{ps88} it is sufficient to show
\begin{equation}\label{040701}
 \lim_{n \rightarrow \infty } \int_B \frac{1}{2} |{\bf y}^n_t({\bf x},t) |^2 + W(\nabla {\bf y}^n({\bf x},t)) d{\bf x} \\
=  \int_B \frac{1}{2} |{\bf y}_t({\bf x},t) |^2 + W(\nabla {\bf y}({\bf x},t)) d{\bf x} + \frac{\omega_d}{d} w(0,t)^d L.
\end{equation}
To this end we decompose the integral on the left hand side of \eqref{040701} into three parts
\begin{multline}
  \int_B \frac{1}{2} |{\bf y}^n_t |^2 + W(\nabla {\bf y}^n) d{\bf x}\\
=  \int_B \frac{1}{2} |{\bf y}^n_t |^2 d{\bf x} + \int_B  W(\nabla {\bf y}^n)  \charf_{\{|{\bf x}| > \frac{1}{n}\}}d{\bf x} 
+ \int_B  W(\nabla {\bf y}^n) \charf_{\{|{\bf x}| < \frac{1}{n}\}}d{\bf x}
=: I^n_1 + I^n_2 +I^n_3.
\end{multline}
The fact that $I^n_1$ converges to $\int_B \frac{1}{2} |{\bf y}_t({\bf x},t) |^2 d{\bf x}$ follows from the boundedness 
of ${\bf y}^n_t$ and its pointwise convergence to ${\bf y}_t$.
Regarding the convergence of $I^n_2$ we note that due to Lemmas \ref{lem:bulkest3d} and \ref{lem:vest3d}  the integrand 
is bounded by $C(1 + \frac{1}{|{\bf x}|} + \frac{1}{|{\bf x}|^ 2}) $ which is integrable.
Moreover, the integrand pointwise converges to $ W(\nabla {\bf y}).$
Thus, it remains to determine the limit of $I_3^n.$
\begin{multline}
 \frac{2}{\omega_d} I^n_3 =  \int_0^{\frac{1}{n}} 2 W(\nabla {\bf y}^n)  R^{d-1} dR\\
=  \int_0^{\frac{1}{n}}\! \lambda_1^n(R,t)^2  R^{d-1} dR +  \int_0^{\frac{1}{n}}\!(d-1) \lambda_2^n(R,t)^2  R^{d-1} dR
+ \int_0^{\frac{1}{n}}\! 2h(v^n(R,t)) R^{d-1} dR =: I^n_{3,1} + I^n_{3,2} +I^n_{3,3}.
\end{multline}
The first two summands converge to zero as the integrand is bounded, quadratic in $n$ and $d \geq 3.$
We rewrite the third summand
\[ \frac{\omega_d}{2} I^n_{3,3} = \omega_d \int_0^1 h(v^n(\frac{R}{n},t)) \frac{R^{d-1}}{n^d} dR\]
such that the integrand can be estimated by 
\[ |  h(v^n(\frac{R}{n},t)) \frac{R^{d-1}}{n^d}| \leq \max\{ h(c), h(c(1 + n^d t^d))\} n^{-d}< C < \infty
\]
independent of $n$ for fixed $t>0.$
Thus, by the dominated convergence theorem we get
\begin{equation}
 \begin{split}
  \lim_{n \rightarrow \infty} \frac{\omega_d}{2}  I^n_{3,3} &= \omega_d \int_0^1 \lim_{n \rightarrow \infty} \frac{h(v^n(\frac{R}{n},t))}{n^d} R^{d-1} dR\\
&= \omega_d \int_0^1 \lim_{n \rightarrow \infty} 
\frac{ h(2^d w(0,t)^d n^d \phi(R) \frac{\Phi(R)^{d-1}}{R^{d-1}} + \mathcal{O}(n^{d-1}))}{\frac{n^d}{R^{d-1}}} dR\\
&=\omega_d \int_0^1 L 2^d w(0,t)^d \phi(R) \Phi(R)^{d-1} dR\\
&= \frac{\omega_d L 2^d w(0,t)^d}{d} [\Phi(R)^d]_0^1 = \frac{\omega_d L  w(0,t)^d}{d}.
 \end{split}
\end{equation}

\qed

Our next aim is to show that, in fact, energy is needed to create the cavity. 
\begin{proposition}\label{lem:increase}
 Let $\lambda>0$ be given, $\lim_{v \rightarrow \infty} \frac{h(v)}{v}< \infty$ and ${\bf y}$ a cavitating solution as computed in \cite{ps88}. Then
\[ \frac{d}{dt} E[{\bf y},B](t) >0\]
for any ball $B$ containing the whole wave fan at time $t.$
\end{proposition}

Before we can prove this  we need a preparatory lemma.
\begin{lemma}
 Under the conditions of Lemma \ref{lem:increase} 
\[ r(\sigma)^d \left(1 - \frac{r'(\sigma-)}{\lambda}\right) - r(0)^d \leq 0,\]
where $\sigma$ is the speed of the shock.
\end{lemma}

\proof 
 This proof heavily relies on the estimates in Lemma  \ref{lem:psest}, in particular the monotonicity of $v$.
Having said that, we calculate
\begin{multline}
  r(\sigma)^d - r(0)^d = \int_0^\sigma \frac{d}{ds} \left( r(s)^d\right) ds = d \int_0^\sigma r'(s) r(s)^{d-1} ds
    = d\int_0^\sigma v(s) s^{d-1}  ds\\
     \leq d\int_0^\sigma r'(\sigma-) \left(\frac{r(\sigma)}{\sigma}\right)^{d-1} s^{d-1}  ds
     = \sigma^d r'(\sigma-)\left(\frac{r(\sigma)}{\sigma}\right)^{d-1} 
    = r'(\sigma-) r(\sigma)^{d-1} \sigma
    = \frac{r'(\sigma-)}{\lambda} r(\sigma)^d 
\end{multline}
as $r(\sigma)=\lambda \sigma.$
\qed

\medskip
\noindent
{\it Proof of Lemma \ref{lem:increase}:}
 By Proposition \ref{lem:energy3d2} we have 
\begin{equation}
 E[{\bf y},B](t) - E[\lambda {\bf x},B](t) = \frac{t^d \omega_d}{d} D
\end{equation}
where 
\begin{multline}
 D:= \sigma^d\Big[\frac{1}{2} w_R(t\sigma-,t)^2 + h(w_R(t\sigma-,t)\lambda^{d-1}) - \frac{1}{2} \lambda^2 - h(\lambda^d)\\
  + \frac{1}{2} \left[   w_R(t \sigma-,t) + h'(w_R(t \sigma-,t) \lambda^{d-1} ) \lambda^{d-1}+ \lambda + h'(\lambda^d)\lambda^{d-1}   
\right]  
(\lambda -  w_R(t \sigma-,t) )\Big] + r(0)^d L.
\end{multline}
Thus, it is sufficient to show $D>0.$
We rewrite $D$ using \eqref{eq:rdef}
\begin{multline}\label{old29}
 D= \sigma^d\Big[\frac{1}{2} r'(\sigma-)^2 + h(r'(\sigma-)\lambda^{d-1}) - \frac{1}{2} \lambda^2 - h(\lambda^d) 
+ \frac{1}{2} (r'(\sigma-) + \lambda)(\lambda -  r'(\sigma-)) \\
  + \frac{\lambda^{d-1}}{2} \left(  h'(\lambda^d) - h'(r'(\sigma-) \lambda^{d-1} )  
\right)  
(\lambda -  r'(\sigma-) )   + h'(r'(\sigma-) \lambda^{d-1} )\lambda^{d-1}(\lambda -  r'(\sigma-) ) \Big] + r(0)^d L.
\end{multline}
We note that
\begin{equation}\label{old29b}
\lambda^{d-1} \left(  h'(\lambda^d) - h'(r'(\sigma-) \lambda^{d-1} ) 
\right)  
(\lambda -  r'(\sigma-) ) \geq 0
\end{equation}
as $a \mapsto h'(a \lambda^{d-1})$ is increasing, because $h$ is convex.
Furthermore, 
\begin{equation}
 \label{old29c}
 h'(r'(\sigma-) \lambda^{d-1} )(\lambda -  r'(\sigma-) ) > 0
\end{equation}
because $\lambda > r'(\sigma-) $ and $0=h'(H)< h'(r'(\sigma-) \lambda^{d-1} ),$ cf. Lemma \ref{lem:psest}.
We combine \eqref{old29}, \eqref{old29b} and \eqref{old29c} to obtain
\begin{equation}\label{old30}
 \begin{split}
  D&> \sigma^d\Big[\frac{1}{2} r'(\sigma-)^2 + h(r'(\sigma-)\lambda^{d-1}) - \frac{1}{2} \lambda^2 - h(\lambda^d) 
+ \frac{1}{2} (r'(\sigma-) + \lambda)(\lambda -  r'(\sigma-)) \Big] + r(0)^d L\\
 &= \sigma^d \left[ h(r'(\sigma-)\lambda^{d-1}) - h(\lambda^d)  \right] + r(0)^d L\\
 &= \sigma^d h'(\tilde v) (r'(\sigma-) - \lambda)\lambda^{d-1} +r(0)^d L
 \end{split}
\end{equation}
 for some $r'(\sigma-) \lambda^{d-1} \leq \tilde v \leq \lambda^d.$
 Using in \eqref{old30} that $L> h'(\tilde v)$, as $h''>0$ and $L=\lim_{v \rightarrow \infty} \frac{h(v)}{v}=\lim_{v \rightarrow \infty} h'(v)$,  we find
 \begin{multline}
   D >  (\sigma^d \lambda^{d-1} (r'(\sigma-) - \lambda) + r(0)^d)L = \left( \sigma^d \lambda^d \left( \frac{r'(\sigma-)}{\lambda}-1\right) +r(0)^d \right) L\\
 = \left( r(0)^d -r(\sigma)^d \left( 1- \frac{r'(\sigma-)}{\lambda} \right) \right) L \geq 0.
 \end{multline}

\qed

\section{Gas dynamics with vacuum in Lagrangean coordinates}
\label{sec-vacuum}
In this section we describe how our notion of slic--solution easily extends to the system of gas dynamics in Lagrangean 
coordinates. In particular, we will
consider a solution to the Riemann problem for the $p$--system containing a vacuum state. The $p$--system is given by
\begin{equation}\label{eq:psys}
 \begin{split}
  u_t - v_x &=0\\
v_t + p(u)_x&=0
 \end{split}
\quad \text{ in } \R \times (0,\infty),
\end{equation}
where $u$ denotes specific volume and $v$ denotes velocity,  for some polytropic gas law, i.e.,
 $ p(u)=\frac{1}{\gamma} u^{-\gamma}$ with $\gamma>1$.
For $\bar u, \bar v >0$ we complement \eqref{eq:psys} with initial data
\begin{equation}\label{eq:vacid}
 u(x,0)= \bar u \text{ for all } x \in \R,\quad v(x,0)=\left\{\begin{array}{ccc}
                                                               -\bar v & \text{for} & x<0\\
                                                                \bar v & \text{for} & x>0
                                                              \end{array}
  \right..
\end{equation}
In Section 9.6 of \cite{dafermosbook} it is noted that whether  \eqref{eq:psys},\eqref{eq:vacid} admits a standard weak solution depends on the sign of 
\[ w:= \bar u^{\frac{1-\gamma}{2}} + \frac{1-\gamma}{2} \bar v.\]
In case $w<0$ there exists no standard weak solution but a solution containing a  vacuum state,
which is given in \cite{dafermosbook} and has the form 
$u = u \Big ( \frac{x}{t} \Big )$, $v = v \Big ( \frac{x}{t} \Big )$ where $\xi = \frac{x}{t}$,  $\xi_F:= \bar u^{-\frac{\gamma +1}{2}}$,
\begin{equation}\label{vac0}
\begin{split}
u(\xi) &= \left\{ \begin{array}{ccc}
                   \bar u   & \text{ for } & \xi \leq -\xi_F\\
                   |\xi|^{-\frac{2}{\gamma +1}} - \frac{4}{\gamma -1} w \delta_0 & \text{ for } & -\xi_F \leq \xi \leq \xi_F\\
                   \bar u   & \text{ for } & \xi_F \leq \xi
                  \end{array}
 \right.\\
v(\xi) &= \left\{ \begin{array}{ccc}
                   -\bar v   & \text{ for } & \xi \leq -\xi_F\\
                   \frac{2}{\gamma-1} \operatorname{sign}(\xi) \left( |\xi|^{\frac{\gamma-1}{\gamma+1}} -w \right) & \text{ for } & -\xi_F \leq \xi \leq \xi_F\\
                   \bar v   & \text{ for } & \xi_F \leq \xi
                  \end{array}
 \right.
\end{split}
\end{equation}
It is noted in  \cite{dafermosbook} that \eqref{vac0} is a distributional solution of \eqref{eq:psys}
provided $u(\xi)^{-\gamma}$ is understood as a continuous function that vanishes at $\xi = 0$, and 
it is referred to an Eulerian description for a validation of this solution.
Moreover, \eqref{vac0} converges pointwise to the initial data \eqref{eq:vacid} for $t \rightarrow 0.$
Subsequently we will show that our approach can treat this solution in a systematic fashion.
To fit \eqref{vac0} into our framework we compute a displacement field $y(x,t)=t r(\frac{x}{t})$ such that
\begin{equation}\label{eq:vac1}
\del_t (t r(\frac{x}{t})) = v(x,t), \quad \del_x (t r(\frac{x}{t})) = u(x,t).
\end{equation}
It is easy to verify that \eqref{eq:vac1} is satisfied for
\begin{equation}\label{eq:vac2}
r(\xi)= \frac{\gamma+1}{\gamma-1} \xi^{\frac{\gamma-1}{\gamma+1}} - \frac{2}{\gamma-1}w \ \text{ for } \xi>0, \quad r(\xi)=-r(-\xi) \ \text{ for } \xi<0.
\end{equation}

\begin{definition}
\label{def:selfsimslic}
Let $\phi$ be any mollifier satisfying the conditions from Definition 2.1 and let $\underset{\xi}{\star}$ denote mollification in $\xi$. 
For  $r \in L^1_{loc} (\R )$ we call $y(x,t)=t \, r \Big ( \frac{x}{t} \Big )$ a self--similar slic--solution of \eqref{eq:psys} if 
 the sequences
\begin{equation}\label{eq:vac3} 
r_n(\xi) := \big (\phi_n \underset{\xi}{\star} r \big )(\xi), \quad u_n(\xi) := r'_n(\xi), \quad v_n(\xi) = r_n(\xi) - \xi r'_n(\xi),
\end{equation}
 fulfill
\begin{equation}
\label{eq:vac8}
\int_{-\infty}^{\infty} (\frac{1}{\gamma}(u_n)^{-\gamma} - \xi v_n)\psi_\xi - v_n \psi d\xi \rightarrow 0 \text{ as } n \rightarrow \infty
\end{equation}
for all mollifiers $\phi$ and all $\psi \in C_c^1(\R)$.
\end{definition}
\begin{remark}
 As $u_n$ and $v_n$ are derived from the same displacement field it is clear that
$(u_n,v_n)$ solves \eqref{eq:psys}$_1$ exactly for every $n \in \mathbb{N}.$
\end{remark}
Our aim is to show that \eqref{vac0} defines a self--similar slic--solution. To do this we need the following estimates.

\begin{lemma}\label{lem:vac:est}
There exists $C>0$ such that
\[ |u_n(\xi)^{-\gamma}|,|v_n(\xi)| < C\quad 
\text{ for all $\xi \in \R$ and $n \in \mathbb{N}.$}\]
\end{lemma}
\proof
The anti--symmetry of $r$ is inherited by  $r_n$ and we may write $r_n$ as follows 
\begin{equation}
\label{eq:vac4}
r_n(\xi) = \int_0^\infty \big( \phi_n(\xi -\tilde \xi) - \phi_n(\xi+\tilde \xi)\big) r(\tilde \xi) d \tilde \xi.
\end{equation}
We will compute expressions for  $u_n(\xi)$ and $v_n(\xi)$ which are well suited for estimates:
\begin{multline}\label{eq:vac5}
 u_n(\xi)= \frac{d}{d\xi} \left( \int_0^\infty \phi_n(\xi-\tilde \xi) ( \frac{\gamma+1}{\gamma-1} \tilde \xi^{\frac{\gamma-1}{\gamma+1}} 
- \frac{2w}{\gamma-1})d\tilde \xi
 - \int_{-\infty}^0 \phi_n(\xi-\tilde \xi)(    \frac{\gamma+1}{\gamma-1} (-\tilde \xi)^{\frac{\gamma-1}{\gamma+1}} 
- \frac{2w}{\gamma-1})d\tilde \xi \right)\\
= - \int_0^\infty \frac{d}{d \tilde \xi}\phi_n(\xi-\tilde \xi) ( \frac{\gamma+1}{\gamma-1} \tilde \xi^{\frac{\gamma-1}{\gamma+1}} 
- \frac{2w}{\gamma-1})d\tilde \xi
 + \int_{-\infty}^0 \frac{d}{d \tilde \xi}\phi_n(\xi-\tilde \xi)(    \frac{\gamma+1}{\gamma-1} (-\tilde \xi)^{\frac{\gamma-1}{\gamma+1}} 
- \frac{2w}{\gamma-1})d\tilde \xi \\
= \int_{-\infty}^\infty  \phi_n(\xi - \tilde \xi) |\tilde \xi|^{-\frac{2}{\gamma+1}} d \tilde \xi - \frac{4w}{\gamma-1} \phi_n(\xi).
\end{multline}
%
Due to \eqref{eq:vac4} we find for $0 < \xi < \xi_F$
\begin{equation}\label{eq:vac6}
\begin{split}
v_n(\xi)&=r_n(\xi) - \xi r'_n(\xi)\\
&=\int_0^\infty \big( \phi_n(\xi -\tilde \xi) - \phi_n(\xi+\tilde \xi)\big) 
\big( \frac{\gamma+1}{\gamma-1} \tilde \xi^{\frac{\gamma-1}{\gamma+1}} - \frac{2}{\gamma-1}w\big) d \tilde \xi
\\
&\qquad - \int_0^\infty \big( \phi_n(\xi -\tilde \xi) + \phi_n(\xi+\tilde \xi)\big) \xi \tilde \xi^{-\frac{2}{\gamma+1}} d\tilde \xi + \frac{4}{\gamma-1} w \phi_n(\xi)\xi\\
&=\int_0^\infty \big( \phi_n(\xi -\tilde \xi) - \phi_n(\xi+\tilde \xi)\big) 
\big( \frac{2}{\gamma-1} \tilde \xi^{\frac{\gamma-1}{\gamma+1}}- \frac{2}{\gamma-1}w\big) d \tilde \xi
\\
&\qquad + \int_0^\infty \big( \phi_n(\xi -\tilde \xi)(\tilde \xi - \xi) - \phi_n(\xi+\tilde \xi)(\xi+\tilde \xi)\big) 
\tilde \xi^{-\frac{2}{\gamma+1}} d\tilde \xi + \frac{4}{\gamma-1} w \phi_n(\xi)\xi.
\end{split}
\end{equation}

\noindent
To show that $\frac{1}{\gamma} u_n^{-\gamma}$ is bounded it is sufficient to show that for $|\xi| < \xi_F$ the specific volume $u_n(\xi)$ 
is bounded from below by some positive constant.
As $\gamma>1$ and $w<0,$ we find for $|\xi| < \xi_F$ that
\begin{multline}
 u_n (\xi) \geq \int_{\xi-\frac{1}{n}}^{\xi+\frac{1}{n}}  \phi_n(\xi - \tilde \xi) |\tilde \xi|^{-\frac{2}{\gamma+1}} d \tilde \xi
\geq \int_{\xi-\frac{1}{n}}^{\xi+\frac{1}{n}}  \phi_n(\xi - \tilde \xi) (|\xi| +1)^{-\frac{2}{\gamma+1}} d \tilde \xi\\
= (|\xi| +1)^{-\frac{2}{\gamma+1}} \geq (\xi_F +1)^{-\frac{2}{\gamma+1}}.
\end{multline}

Regarding $v_n(\xi)$ we find for $\xi \in [-\xi_F, \xi_F]$
\begin{align*}
|v_n(\xi)| &\leq    \Big( \frac{2}{\gamma-1} ( \xi_F +1)^{\frac{\gamma-1}{\gamma+1}}\Big)
 \int_0^\infty \big( \phi_n(\xi -\tilde \xi) + \phi_n(\xi+\tilde \xi)\big) 
 d \tilde \xi
- \frac{2}{\gamma-1}w\\
 &\qquad + \int_0^\infty 2 |\phi|_\infty \charf_{\{\xi -\frac{1}{n} < \tilde \xi < \xi +\frac{1}{n} \}}
 \tilde \xi^{-\frac{2}{\gamma+1}} d\tilde \xi - \frac{4}{\gamma-1} w 
 |\phi|_\infty\\
 &\leq \frac{2}{\gamma -1}( \xi_F +1)^{\frac{\gamma-1}{\gamma+1}}
 - \frac{2}{\gamma-1}w + 2 |\phi|_\infty \int_0^{ \xi_F +1 } \tilde \xi^{-\frac{2}{\gamma+1}} d \tilde \xi
 - \frac{4}{\gamma-1} w |\phi|_\infty\\
 &= \frac{2}{\gamma-1} \big( 1 + \frac{\gamma-1}{\gamma+1} |\phi|_\infty\big) (\xi_F +1)^{\frac{\gamma-1}{\gamma+1}}
 - \frac{2}{\gamma-1}w- \frac{4}{\gamma-1} w |\phi|_\infty < \infty.
\end{align*}
\qed

\begin{theorem}\label{lem:vac:sol}
 The function $y(x,t):=t r(\frac{x}{t})$ with $r$ defined by  \eqref{eq:vac2} is a self--similar slic--solution of \eqref{eq:psys}.
\end{theorem}
\proof
As the constant state and the rarefaction wave in $(u,v)$ given in \eqref{vac0} are linked continuously it is sufficient to consider test functions $\psi$ 
with support inside the rarefaction/vacuum wave, i.e.,
$\supp(\psi) \subset [-\xi_F,\xi_F]$.
Furthermore, observe that
$$
\begin{aligned}
(u_n)^{-\gamma} (\xi)  &\to  |\xi|^{\frac{2\gamma}{\gamma +1}} \, ,
\\
v_n(\xi)  &\to 
 \frac{2}{\gamma-1} \operatorname{sign}(\xi) \left( |\xi|^{\frac{\gamma-1}{\gamma+1}} -w \right) \, ,
\end{aligned}
$$
pointwise. Using lemma \ref{lem:vac:est} and the dominated convergence theorem, we obtain
\begin{align*}
&\lim_{n \rightarrow \infty} \int_{-\infty}^{\infty} (\frac{1}{\gamma}(u_n)^{-\gamma} - \xi v_n)\psi_\xi - v_n \psi d\xi \\
&=\int_{-\infty}^{\infty}\lim_{n \rightarrow \infty} ( (\frac{1}{\gamma}(u_n)^{-\gamma} - \xi v_n)\psi_\xi - v_n \psi) d\xi \\
&=\int_{0}^{\infty} \left( \frac{1}{\gamma} (\xi^{\frac{-2}{\gamma+1}})^{-\gamma} - \xi \frac{2}{\gamma-1} \xi^{\frac{\gamma-1}{\gamma+1}} 
+ \xi \frac{2}{\gamma-1}w \right) \psi_\xi(\xi)
- \left(  \frac{2}{\gamma-1} \xi^{\frac{\gamma-1}{\gamma+1}} - \frac{2}{\gamma-1}w\right) \psi(\xi) d\xi\\
&+ \int_{-\infty}^0 \!\! \left( \frac{1}{\gamma} ((-\xi)^{\frac{-2}{\gamma+1}})^{-\gamma} + \xi \frac{2}{\gamma-1} (-\xi)^{\frac{\gamma-1}{\gamma+1}} 
- \xi \frac{2}{\gamma-1}w \right) \psi_\xi(\xi)
+ \left(  \frac{2}{\gamma-1} (-\xi)^{\frac{\gamma-1}{\gamma+1}} - \frac{2}{\gamma-1}w\right) \psi(\xi) d\xi  \\
&=- \int_{0}^{\infty} \left( \frac{2}{\gamma+1} \xi^{\frac{\gamma-1}{\gamma+1}} - \frac{2}{\gamma-1} \xi^{\frac{\gamma-1}{\gamma+1}}
- \xi \frac{2}{\gamma+1} \xi^{\frac{-2}{\gamma+1}} + \frac{2w}{\gamma-1} + \frac{2}{\gamma-1} \xi^{\frac{\gamma-1}{\gamma+1}} 
- \frac{2w}{\gamma-1}\right)\psi(\xi) d\xi\\
&-\! \int_{-\infty}^{0}\!\! \left(\! \frac{-2}{\gamma+1} (-\xi)^{\frac{\gamma-1}{\gamma+1}} + \frac{2}{\gamma-1}(- \xi)^{\frac{\gamma-1}{\gamma+1}}
+ \frac{2}{\gamma+1} (-\xi)^{\frac{\gamma-1}{\gamma+1}} - \frac{2w}{\gamma-1} - \frac{2}{\gamma-1} (-\xi)^{\frac{\gamma-1}{\gamma+1}} 
+ \frac{2w}{\gamma-1}\!\right)\!\psi(\xi) d\xi\\
&=0.
\end{align*}
\qed

\noindent We define the energy of a self--similar slic--solution inside a wave fan in a manner  analogous 
to Definition \ref{def:energy3d}.
\begin{definition}
For any $\bar \xi>0$ the energy of a self--similar slic--solution $y$ inside the wave fan from $-\bar \xi$ to $\bar \xi$ is defined by
 \[ E[ y  \, ,  \, (-\bar \xi ,\bar \xi)  ]:= \lim_{n \rightarrow \infty}\int_{-\bar \xi}^{\bar \xi} W(u_n) + \frac{1}{2} (v_n)^2 d\xi \quad 
\text{ with } W(u) = \frac{1}{\gamma(\gamma-1)} u^{1-\gamma}.\] 
\end{definition}

\begin{proposition}\label{lem:vac:energy}
  There is no contribution of the vacuum state to the energy. In particular,
\[  E[ y \, ,  \, (-\bar \xi ,\bar \xi) ]=\int_{-\bar \xi}^{\bar \xi} \frac{1}{\gamma(\gamma -1)} |\xi|^{\frac{2(\gamma-1)}{\gamma+1}} +
  \frac{2}{(\gamma-1)^2} \left( |\xi|^{\frac{\gamma-1}{\gamma+1}} -w\right)^2 d\xi
\]
for any $0 <\bar \xi< \xi_F$.
\end{proposition}

\proof
Due to Lemma \ref{lem:vac:est}
 the integrand $W(u_n) + \frac{1}{2} (v_n)^2$ is uniformly bounded and therefore limit and integral commute such that
\begin{multline}
\lim_{n \rightarrow \infty}\int_{-\bar \xi}^{\bar \xi} W(u_n) + \frac{1}{2} (v_n)^2 d\xi =
\int_{-\bar \xi}^{\bar \xi} \lim_{n \rightarrow \infty} \big(W(u_n) + \frac{1}{2} (v_n)^2\big) d\xi \\
=\int_{-\bar \xi}^{\bar \xi} \frac{1}{\gamma(\gamma -1)} |\xi|^{\frac{2(\gamma-1)}{\gamma+1}} + \frac{2}{(\gamma-1)^2} \left( |\xi|^{\frac{\gamma-1}{\gamma+1}} -w\right)^2 d\xi.
\end{multline}
\qed

\noindent
{\bf Acknowledgements} 
This research was supported by the EU FP7-REGPOT project "Archimedes Center for
Modeling, Analysis and Computation". AET is partially supported by the "Aristeia" program of the Greek
Secretariat for Research.


\begin{thebibliography}{10}




\bibitem{ball82}
{\sc J. M. Ball},
Discontinuous equilibrium solutions and cavitation in
              nonlinear elasticity,
   {\em Philos. Trans. Roy. Soc. London Ser. A}, {\bf 306}, (1982)  557--611.


 \bibitem{ciarlet93}
 {\sc P.G. Ciarlet},
 {\em Mathematical Elasticity, Vol. I, Three-dimensional elasticity},
 North Holland, Amsterdam, 1988.


\bibitem{dafermosbook}
{\sc C. Dafermos},
{\em Hyperbolic Conservation Laws in Continuum Physics},
Third Edition.
Grundlehren der Mathematischen Wissenschaften, 325.
Springer Verlag, Berlin, 2010.

\bibitem{E00}
{\sc G. Ercole},
Delta-shock waves as self-similar viscosity limits. 
{\it Quart. Appl. Math.} {\bf 58} (2000), 177 - 199.

\bibitem{DS05}
{\sc V.G. Danilov and V.M. Shelkovich},
Delta-shock wave type solution of hyperbolic systems of conservation laws.
{\em Quart. Appl. Math.} {\bf 63} (2005), 401 - 427.

\bibitem{KK89}
{\sc B.L. Keyfitz and H.C. Kranzer},
A viscosity approximation to a system of conservation laws with no classical
Riemann solution. 
{\it Lecture Notes in Math}. No. {\bf 1402},  Berlin: Springer, 1989,  pp. 185 - 197. 

\bibitem{L09}
{\sc O. Lopez-Pamies}
Onset of Cavitation in Compressible, Isotropic, Hyperelastic Solids
{\em J. Elasticity} {\bf 94}  (2009), 115 - 145.


\bibitem{NS11}
{\sc P.V. Neg{\'r}on Marrero and J. Sivaloganathan},
The radial volume derivative and the critical boundary displacement for cavitation. 
{\em SIAM J. Appl. Math.} {\bf 71} (2011), 2185 - 2204. 
 
\bibitem{ps88}
{\sc K.A. Pericak-Spector and S.J. Spector},
Nonuniqueness for a hyperbolic system: cavitation in nonlinear elastodynamics.
{\em Arch. Rational Mech. Anal.} {\bf 101} (1988), 293 - 317.

\bibitem{ps98}
{\sc K.A. Pericak-Spector and S.J. Spector},
Dynamic cavitation with shocks in nonlinear elasticity.
{\em Proc. Royal Soc. Edinburgh Sect A} {\bf 127} (1997), 837 - 857.

\bibitem{SS03}
{\sc J. Sivaloganathan and S.J. Spector},
Myriad radial cavitating equilibria in nonlinear elasticity.
{\em SIAM J. Appl. Math.} {\bf 63} (2003), 1461 - 1473.

\bibitem{TN65}
C. Truesdell, W. Noll, {\it The Non-Linear Field Theories of
Mechanics}. Handbuch der Physik, III$/3$, Berlin: Springer, 1965.


   
  
\end{thebibliography}
\end{document}